\documentclass[12pt]{amsart}

\usepackage{amsfonts,amsmath,latexsym,amssymb,verbatim,amsbsy,amsthm}

\usepackage{graphicx}
\usepackage{nicefrac}
\usepackage{dsfont}
\usepackage{mathrsfs}
\usepackage{cancel}
\usepackage[normalem]{ulem}
\usepackage{color}
\usepackage{pifont}

\usepackage{array}
\usepackage{slashbox}

\usepackage{setspace} %set stretch
\usepackage{caption,subcaption} %figure's stuff

\usepackage{algorithm} %algorithm

\usepackage{algpseudocode} %algorithmic

\usepackage{grffile}  %align figures
\usepackage[top=1in, bottom=1in, left=1in, right=1in]{geometry}

\usepackage[dvipsnames]{xcolor}
\usepackage{relsize}

\usepackage[colorlinks=true, pdfstartview=FitV, linkcolor=RoyalBlue,citecolor=ForestGreen, urlcolor=blue]{hyperref}

\numberwithin{equation}{section}
\numberwithin{figure}{section}
\numberwithin{table}{section}
\usepackage{float}

\renewcommand{\geq}{\geqslant}

\renewcommand{\leq}{\leqslant}

\newcommand{\ds}{\displaystyle} 
\newcommand{\be}{\begin{equation}}
\newcommand{\ee}{\end{equation}}

\theoremstyle{plain}
\newtheorem{THEOREM}{Theorem}[section]

\newtheorem{theorem}[THEOREM]{Theorem}

\newtheorem{proposition}[THEOREM]{Proposition}

\theoremstyle{definition}

\theoremstyle{remark}

\newtheorem{remark}[THEOREM]{Remark}

\newcommand{\mymax}[1]{{}^{#1}_{+}} % instead of \max

\newcommand{\imin}{{i_n}}

\newcommand{\relmin}{\widetilde{m}^{n+1}_i} 
\newcommand{\relminin}{\widetilde{m}^{n+1}_{i_n}} 
\newcommand{\R}{\mathbb{R}}

\newcommand{\calN}{{\mathcal N}}

\newcommand{\bX}{\mathbf{X}}
\newcommand{\bx}{\mathbf{x}}
\newcommand{\bxn}{\bx^n}

\newcommand{\bxmin}{\bx^*}
\newcommand{\bXnmin}{\bX^n_-}
\newcommand{\bXn}{\bX^n}
\newcommand{\bXnplusmin}{\bX^{n+1}_-}
\newcommand{\bXnmina}{\bX^{n_\alpha}_-}
\newcommand{\bXna}{\bX^{n_\alpha}}

\newcommand{\bxin}{\bx_i^n}

\newcommand{\hin}{h^n_i}

\newcommand{\etain}{\eta^n_i}
\newcommand{\mq}{q}
\newcommand{\etaintwo}{\etain} 
\newcommand{\bp}{{\mathbf p}}

\newcommand\bpin{\bp^n_i}
\newcommand\bvin{{\mathbf v}^n_i}
\newcommand\bqin{{\mathbf q}^n_i}
\newcommand{\qin}{q^n_i}
\newcommand\boin{{\boldsymbol \omega}^n_i}

\newcommand{\by}{\mathbf{y}}

\newcommand{\bz}{\mathbf{z}}

\newcommand{\bv}{\mathbf{v}}

\newcommand{\bY}{\mathbf{Y}}
\newcommand{\hf}{\frac{1}{2}}
\newcommand{\algo}{Algorithm }

\newcommand{\gam}{\gamma}

\newcommand{\clam}{\lambda} %instead of \clam
\newcommand{\cdel}{\delta}

\def\SBD{SBRD }

\DeclareMathOperator*{\argmin}{argmin}
\DeclareMathOperator*{\argmax}{argmax}

\begin{document}
\title{Swarm-Based Optimization with Random Descent}

\dedicatory{Dedicated to Shi Jin and for many years of friendship}

\author{Eitan Tadmor}
\address{Department of Mathematics and Institute for Physical Science \& Technology\newline \hspace*{0.3cm}   University of Maryland, College Park}
\email{tadmor@umd.edu\newline
\hspace*{0.3cm} ORCID: 0000-0001-7424-6327}

\author{Anil Zenginoglu}
\address{Institute for Physical Science \& Technology, 
 University of  Maryland, College Park}
\email{anil@umd.edu}

\date{\today}



\thanks{\textbf{Acknowledgment.} Research was supported in part by ONR grant N00014-2112773.}

\begin{abstract}
We extend our study  of the swarm-based gradient descent method for non-convex optimization, \cite{LTZ2022}, to allow random descent directions. We recall that the swarm-based approach consists of a swarm of agents, each identified with a position, $\bx$, and mass, $m$. The key is the transfer of mass from high ground to low(-est) ground. The mass of an agent dictates its step size: lighter agents take larger steps. In this paper, the essential new feature is the choice of direction: rather than restricting the swarm to march in the steepest gradient descent, we let agents proceed in randomly chosen directions centered around --- but otherwise different from --- the gradient direction. The random search secures the descent property while at the same time, enabling greater exploration of ambient space. Convergence analysis and benchmark optimizations demonstrate the effectiveness of the swarm-based random descent method as a multi-dimensional global optimizer.
\end{abstract}

\maketitle
\setcounter{tocdepth}{1}
\tableofcontents

%\ifx%%%
\section{Introduction. The importance of random marching directions.}

\noindent
In this work we extend  our study of swarm-based approach for non-convex optimization, \cite{LTZ2022}, with the aim of finding minimizer(s) of a loss function, $\argmin_{\bx\in \Omega}F(\bx)$, over an ambient bounded  set $\Omega \subset \R^d$. Classical iterative algorithms for numerical optimization  employ a single agent  which   explores the ambient space by  successively improving the position of approximate optimize(s), e.g.,  \cite{nocedal1999numerical,boyd2004convex,nocedal2006conjugate} and the references therein or the more recent \cite{kingma2017adam,
liu2022adaptive} etc. Unlike those single-agent iterations,   the swarm-based methods use a \emph{crowd of coordinated agents} --- the swarm,  to explore $\Omega$,  e.g., \cite{SA1,SA2,PSO,CBO1,CBO-analytical,CBO-timediscrete,carrillo2021consensus,totzeck2022trends,LTZ2022,grassi2023pso}. Here we follow the  swarm-based approach introduced in \cite{LTZ2022}, in which  the swarm consists of $N$ agents, each is identified with a  position $\bxin$, \underline{and} an independent mass (or weight), $m^n_i$, 
 \[
 \bxin=\bx_i(t^n) \in \Omega \subset \R^d, \ \  m_i^n=m_i(t^n)\in (0,1], \qquad i=1,2,\ldots, N.
 \] 
Thus, the distinctive feature of our swarm-based iterations is the fact that they are  embedded in the larger space, $(\bxin,m^n_i)\in \Omega\times [0,1]\subset \R^{d+1}$: the additional mass-parametrization  is the essential platform  which  enables  proper  coordination of agents, in order to improve the overall configuration of the swarm  in its search for an optimizer.\newline
The basic step takes the form
\[
\bx_i^{n+1}= \bxin-\hin\bpin, \qquad i=1,2,\ldots, N.
\] 
It reflects the move of an agent  from its current position, $\bxin$, in direction $\bpin$ with step-size $\hin$. In \cite{LTZ2022} we advocated the use of gradient direction, $\bpin=\nabla F(\bxin)$. By properly adjusting the step-size  which takes into account  the weights  of all other agents in the crowd, $\hin=h_i\big(\bxin, \{m^n_j\}_{1\leq j\leq N}\big)$ (this is where the communication between agents enters),  one is able to secure the all important \emph{descent property}, \cite[eq. (5.5)]{LTZ2022}
\[
F(\bxin-\hin\bpin)\leq F(\bxin)-\clam_i \hin|\nabla F(\bxin)|^2, \qquad \bpin=\nabla F(\bxin);
\]
here $\clam_i$ are the descent amplitudes, depending on the
weights $ \clam_i=\clam_i\big(\{m^n_j\}_{1\leq j\leq N}\big) \in (0,1)$. In this work, we abandon the use of the gradient direction, $\bpin=\nabla F(\bxin)$,  and instead focus on the descent property as the sole guidance in our swarm-based iterations. This allows us to choose from a large cone of directions which still secure the descent property.
By randomly choosing  proper directions $\bpin$, which are still compatible with the descent property, we significantly increase the heterogeneity of the swarm-based method in exploring larger portion of the ambient space, while keeping the overall decent property.\newline
Our  swarm-based method end up with an interplay between positions and weights with takes the schematic description
\[
\begin{split}
\big\{\{m^n_j\}_{1\leq j\leq N}, F(\bxin)\big\} & \mapsto \ m_i^{n+1},\\
\big\{m^{n+1}_i, \bxin\big\} & \mapsto \ \bx_i^{n+1}.
\end{split}
\]
The method  repeatedly transfers mass from high to lower ground  while on the way, driving agents to smaller (lower) loss values; in particular, $\{\min_i F(\bxin)\}$ forms a decreasing  sequence in time, 
which ideally approaches the  global minimizer in the region explored by these agents,
\[
\bxin \stackrel{n\rightarrow \infty}{\longrightarrow}\argmin_{\bx\in \Omega} F(\bx).
\]
The last statement applies to certain sub-sequence, 
$\{\bx^n_{i_n}\}_{n>n_0}$, which is made precise in the main convergence results of theorems \ref{thm:SBD convergence} and \ref{thm:GD-backtracking with PL} below.\newline
A detailed description of this two-stage swarm-based mechanism now follows.

\smallskip\noindent
{\bf Mass transfer}. In the first stage, positions change the distribution of mass:  each agent with mass $m^n_i$ transfers a fraction of its mass, $\etaintwo m^n_i$, to the current global minimizer positioned at $\bx_\imin$ where $\imin:=\argmin_i F(\bxin)$.
\begin{equation}\label{eq:etai}
\left\{\begin{array}{lll}
         \ m_i^{n+1} & =m^n_i  -\etaintwo m^n_i, & i\neq \imin \\ \\
          \ m_\imin^{n+1} & =\displaystyle m_\imin^n +\sum \limits_{i\neq \imin}\etaintwo m_i^{n}, & 
     \end{array} \right. \quad \etain := \Big(\frac{F(\bxin)-F^n_{\textnormal{min}}}{F^{n}_{\textnormal{max}}-F^{n}_{\textnormal{min}}}\Big)^\mq \in (0,1].
\end{equation}
The fraction of mass transfer,  $\Big(\frac{F(\bxin)-F^n_{\textnormal{min}}}{F^{n}_{\textnormal{max}}-F^{n}_{\textnormal{min}}}\Big)^\mq$,  is determined by the \emph{relative height} of each agent, relative to the global extremes\footnote{To prevent vanishing denominator in the extreme case $F_{\textnormal{max}} = F_{\textnormal{min}}$, we adjust \eqref{eq:etai} with a small $\epsilon$-correction, 
 $\displaystyle \etain := \Big(\frac{F(\bxin)-F^n_{\textnormal{min}}}{F^n_{\textnormal{max}}-F^n_{\textnormal{min}}+\epsilon}\Big)^\mq$.
 }, $\displaystyle F^n_{\textnormal{min}}= \min_j F(\bx^n_j)$ and $\displaystyle F^n_{\textnormal{max}}= \max_j F(\bx^n_j)$, and depending on user-choice of a mass transfer parameter,  $\mq\geq 1$. The higher $\mq$ is, the more tamed is the transfer of mass. A systematic study  reported in section \ref{sec:high-q} below reveals   a dramatic improvement when increasing the mass transfer parameter $\mq=2,4, 8$.\newline
Observe that while the total mass is conserved, say $\sum_i m_i^n=1$,  individual masses are redistributed from high to lower ground --- the higher the agent, the larger fraction of its mass will be lost in favor of the agent at the lowest ground. In fact,  the  highest agent in each iteration is eliminated: to be precise, the worst performing agents are eliminated whenever   $1-\eta^n_i = {\mathcal O}(\epsilon) \ll 1$. 
This  is an aggressive `survival of the fittest' protocol,  so that after $N$ iterations the swarm consists of a single agent which should be in the best position to approach the minimum of the space explored so far by the swarm. We note in passing that one   can  adopt a more flexible protocol which allows the worst (highest) agents to survive a few iterations before  elimination; this flexibility would improve the overall success rates of the swarm at the expense of efficiency.\newline
In particular, the dynamic adjustment of masses  in \eqref{eq:etai}  can be interpreted as a particular case of \emph{alignment dynamics}, with `aggressive' protocol in which agents steer towards the   \emph{minimal} heading, $m^n_-=\min_i m^n_i$. Instead  one may consider a more tamed alignment towards an average heading, as originated in \cite{reynolds1987flocks}, see e.g., \cite{tadmor2021mathematics} and the references therein. In this context we refer to the stochastic-based Consensus Based Optimization, \cite{CBO1,CBO-analytical,  CBO-semidiscrete, totzeck2022trends},  steering towards  a properly weighted convex combination, 
$\bar{m}^n=\sum_j \theta^n_j \bx^n_j$ with weights $\theta^n_j= exp(-\alpha F(\bx_j^n))/\big(\sum_k exp(-\alpha F(\bx_k^n))\big)$,  which in turn is driven to a global minimum by letting $\alpha \gg1$. We note that a main novelty  in our approach  is the use of masses $\{m^n_i\}$ in \eqref{eq:etai} as \emph{independent variables}, which evolve alongside the dynamics of positions $\{\bxin\}$ outlined in \eqref{eq:SBD},\eqref{eqs:orient-omega} below. Indeed, carrying out the perspective of alignment dynamics,  the masses provide the essential added dimension as a platform for communication among the agents of the swarm. 

\smallskip\noindent
{\bf Stepping in descent direction --- a random choice approach}. In the second stage, the distribution of mass affects the change of positions,
 \begin{equation}\label{eq:SBD}
 \bx_i^{n+1}= \bxin-\hin\bpin.
 \end{equation}
 %\end{subequations}
 The driving force behind the protocol for choosing the direction, $\bp=\bpin$, and the step size, $h=\hin$, is to secure the following descent property, depending on the relative mass $\relmin$ and a descent parameter $\clam<1$,  
\begin{equation}\label{eq:descentclam}
 F(\bxin-h\bpin)\leq F(\bxin)-\hf\clam \relmin h |\nabla F(\bxin)|^2, \qquad 
 \relmin=\frac{m_i^{n+1}}{\max_i m_i^{n+1}}, \quad \clam<1.
\end{equation}
We recall  that the choice of the gradient direction, $\bpin=\nabla F(\bxin)$,   secures the sharper steepest descent property,  
\[
F(\bxin-h\bpin)_{|\bpin=\nabla F(\bxin)}\leq F(\bxin)-\clam \relmin h|\nabla F(\bxin)|^2,
\]
 which was the basis for the swarm-based gradient descent (SBGD) method we introduced in \cite{LTZ2022}. 
The purpose of this work is to extend the SBGD method by allowing a larger set of descent directions: the emphasis is no longer on the steepest descent along the gradient direction but instead,  allowing a more effective exploration of the ambient space using a \emph{random choice of directions}, $\{\bpin\}$, that still maintains (half the steepest) descent property. This implies that the swarm, stepping in  other than the  gradient direction, will explore  a larger portion  of the ambient space,  which in turn leads to a more effective  search, and proved to be particularly relevant in high-dimensional optimizations, see the numerical simulations reported in section \ref{sec:results}. We refer to this new version as the Swarm-Based Random Descent (SBRD) method.

\vspace*{-0.1cm}
\begin{figure}[ht]
    \centering
    \includegraphics[angle=345,origin=c,scale = 1.2]
    {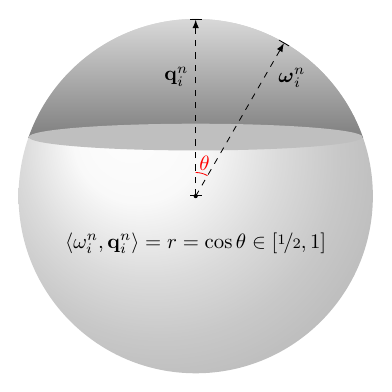}
    \vspace*{-0.7cm}
       \caption{$\bqin$  is the gradient orientation --- the unit vector along the gradient direction, $\nabla F(\bxin)$,  and  the unit vector, $\boin$, is determined by a randomly chosen point on a spherical cap centered around $\bqin$ (shown as the shaded part of the sphere).} \label{fig:spherical_cap}
\end{figure}

\noindent
We provide below a detailed, self-contained description of the SBRD. The stepping protocol is based on  the  choice of direction and step-size. The direction $\bpin$, is determined by its \emph{orientation}, $\boin$,
\begin{subequations}\label{eqs:orient-omega}
\begin{equation}
\bpin=|\nabla F(\bxin)|\boin, \qquad \boin \in {\mathbb S}^{d-1},
\end{equation}
relative to the orientation of the gradient, $\ds \bqin=\frac{\nabla F(\bxin)}{|\nabla F(\bxin)|} \in {\mathbb S}^{d-1}$,
so that
\begin{equation}\label{eq:orient-omega}
\left\langle \boin,\bqin\right\rangle =r,  \qquad \bqin:=\frac{\nabla F(\bxin)}{|\nabla F(\bxin)|}.
\end{equation}
Here, $r$ is  randomly chosen  number from a uniform distribution    in an interval dictated by the relative mass, 
\be
r\in {\mathcal U}\Big[\nicefrac{1}{2}(1+\relmin), 1\Big].
\ee
\end{subequations}
This means that the orientation of $\bpin$ lies in a spherical cap centered around $\bqin$.  The `opening' of the corresponding  spherical cone,   see Figure \ref {fig:spherical_cap}, $\theta:=\arccos \big(\frac{1}{2}(1+\relmin)\big)$. It  is larger for lighter agents, and coincides with the gradient direction, $\nabla F(\bxin)$, for the heaviest agent where $\relmin=1$. The protocol for randomly selecting $\boin$ subject to \eqref{eq:orient-omega} is outlined in section \ref{sec:random} below.

\smallskip\noindent
{\bf Choosing the step size --- a backtracking protocol}. 
It follows that the new position, $\bx^{n+1}(h)=\bxin-h\bpin$ --- viewed as a function of the step size $h$, satisfies the desired descent property, at least for small enough $h$. Indeed,    \eqref{eqs:orient-omega} implies 
  \begin{equation}\label{eq:porient}
  \langle \bpin, \nabla F(\bxin)\rangle= r|\nabla F(\bxin)|^2 \geq \hf(1+\relmin)|\nabla F(\bxin)|^2.
  \end{equation}
 Hence, if we let $L:=\max_{\bx \in \Omega}\|D^2 F(\bx)\|$ with $L< \infty$ serves as a Lipschitz bound\footnote{In fact, one can use a  Lipscitz bound localized to the neighborhood of $\bxin$ that is being visited by the \SBD iterations, but since this neighborhood is not quantified we address the global Lipschitz bound $\ds L=\max_{\alpha,\beta, \bx\in \Omega}\Big|\frac{\partial^2 F(\bx)}{\partial x_\alpha\partial x_\beta}\Big|$.} of $\nabla F$,  then for every $\clam<1$ and  $h<\nicefrac{1}{L}$,  there holds
  \begin{equation}\label{eq:halfdescent}
  \begin{split}
  F(\bx_i^{n+1}(h)) &\leq  F(\bxin)-h\langle \bpin,\nabla F(\bxin)\rangle +
  \frac{h^2}{2}L|\nabla F(\bxin)|^2 \\
   & \leq F(\bxin)-\frac{1}{2}\big(1+\relmin-Lh\big)h|\nabla F(\bxin)|^2\\
   & < F(\bxin)-\frac{1}{2} \clam\relmin h|\nabla F(\bxin)|^2.
   \end{split}
  \end{equation}
Thus, we recover \eqref{eq:descentclam} for any $h<\nicefrac{1}{L}$.\footnote{In fact,  a slightly larger threshold, $h<\frac{1+\relmin(1-\clam)}{L}$, still allows \eqref{eq:halfdescent} to hold.} In particular step size of order $\lesssim \nicefrac{1}{L}$ need not necessary be very small to enforce the descent property. 
 We note, however, that since we have no access to the Lipschitz  bound $L$, we therefore do not have an effective protocol for computing a  step size which secures \eqref{eq:halfdescent}, beyond making the generic statement that it holds for  `sufficiently small' $h$. In fact, $h<\nicefrac{1}{L}$ need not be small and we are  interested in  a  protocol  that identifies the \emph{largest}  $h$ for which \eqref{eq:halfdescent} holds (observe that  the larger $h$ is, the larger is the  descent bound quoted in \eqref{eq:halfdescent}). To this end we use a \emph{backtracking protocol} outlined in \algo \ref{alg:backtracking} below. 
The backtracking algorithm  produces a time step, $h=\hin$,  depending on the position of the agent $\bxin$, and its relative mass $\relmin$,
\[
\hin=h(\bxin,\clam\relmin).
\]
It  secures  the \emph{lower} bound $\displaystyle \hin\geq \frac{\gamma}{L}$ for some $\gamma<1$, so that using \eqref{eq:halfdescent} we finally end up with the descent property of the form,
 \[
F(\bx_i^{n+1}) \leq F(\bxin)-\frac{\gamma}{2L}\clam\relmin|\nabla F(\bxin)|^2, \qquad \bx_i^{n+1}=\bxin-\hin\bpin.
\]
\begin{remark} 
This should be  compared with the descent property of SBGD method restricted to the steepest descent  direction  $\bpin=\nabla F(\bxin)$, for which   we have, \cite[Proposition 5.2]{LTZ2022},
\[
F(\bx_i^{n+1}) \leq F(\bxin)-\frac{2\gamma}{L}\big(1-\clam\relmin\big)\clam\relmin|\nabla F(\bxin)|^2.
\]
Thus, our stepping protocol retains at least half of the steepest descent, while gaining greater heterogeneity in space exploration. In particular, while heavier agents are still restrained by smaller time steps, lighter agents are now allowed to  take larger time steps from a richer set of directions which are aligned with --- but otherwise different from, the gradient  direction. This `greedy' exploration of the ambient space by lighter agents, increases their likelihood of encountering a new neighborhood with a better minimum, which may place one of them as the new heaviest minimizer and so on. 
\end{remark}

\subsection{Why randomization is important} We compare the swarm-based  method 
\begin{equation}\label{eq:swarm-based}
\bx_i^{n+1} = \bxin-h(\bxin,\clam\relmin) \bpin, \qquad i=1,2,\ldots, N,
\end{equation}
in two scenarios: with the gradient direction for SBGD, $\bpin=\nabla F(\bxin)$ and with the randomized direction for SBRD, $\bpin=|\nabla F(\bxin)|\boin$ in \eqref{eqs:orient-omega}. The same backtracking protocol was implemented in both cases. The advantage of randomization in exploring larger regions becomes apparent in \SBD when the number of agents is larger than the dimension of the search space, $N>d$. The results recorded in Table \ref{tab:SBD multiD ackley} for the Ackley function show that \SBD optimization outperforms SBGD optimization in higher dimensions.

\begin{table}[ht]
\setstretch{1.5}
       \centering
    \begin{tabular}{|c| |cc|cc|cc|cc|cc|}
    \hline
    \backslashbox{$d$}{$N$} & \multicolumn{2}{c|}{10} & \multicolumn{2}{c|}{25} & \multicolumn{2}{c|}{50} & \multicolumn{2}{c|}{100} \\
\hline
& \SBD & SBGD & \SBD & SBGD & \SBD & SBGD & \SBD & SBGD\\
\hline
\hline
12 & 13.7\% & 26.7\% & 55.5\%&96.2\% & 88.3\%&100.0\% & 99.2\%&100.0\%\\
13 & 8.8\% & 9.2\% & 49.9\% & 65.5\% & 82.1\% & 95.6\% & 98.1\% & 99.9\% \\
14 & \textbf{3.0\%} & 1.7\% & \textbf{42.4\%} & 22.3\% & \textbf{77.9\%} & 51.0\% & \textbf{96.1\%} & 85.4\% \\
15 & 1.3\% & 0.4\% & \textbf{35.9\%} & 2.7\% & \textbf{70.2\%} & 10.6\% & \textbf{90.5\%} & 23.7\% \\
16 &  0.3\%&0.0\% & \textbf{23.6\%} & 0.1\% & \textbf{60.6\%}&0.8\% & \textbf{85.2\%} & 2.2\%\\
17 & 0.1\% & 0.0\% & \textbf{14.1\%} & 0.0\% & \textbf{50.8\%} & 0.1\% & \textbf{79.1\%} & 0.4\% \\
18 & 0.0\% & 0.0\% & \textbf{8.8\%} & 0.0\% & \textbf{37.3\%} & 0.0\% & \textbf{65.5\%} & 0.0\% \\
19 & 0.0\% & 0.0\% & \textbf{2.0\%} & 0.0\% & \textbf{16.8\%} & 0.0\% & \textbf{48.2\%} & 0.0\% \\
20  & 0.0\%&0.0\% & 0.7\%&0.0\% & \textbf{5.1\%}&0.0\% & \textbf{21.3\%}&0.0\%\\
       \hline
        \end{tabular}
         \smallskip
 \caption{Success rates of \SBD vs. SBGD for global optimization of the $d$-dimensional Ackley function using $N$ agents based on $m=1000$ runs of uniformly generated initial data, $\bx_i^0\in[-3,3]^d$. Backtracking parameters are $\clam=0.2$ and $\gamma=0.9$ (see algorithm \ref{alg:backtracking}). Boldfaced numbers emphasize the cases where \SBD outperforms SBGD by more than 1\%. The randomization provided by \SBD becomes essential beyond the critical dimension $d=13$.}
 \label{tab:SBD multiD ackley}
 \end{table}
\noindent
More can be found in numerical simulations of several benchmark problems presented in section \ref{sec:results}.  

%%%%%%%%%%%%%%%%
\section{Swarm-Based Random Descent (SBRD).  Implementation of   algorithm}\label{sec:construct}
%%%%%%%%%%%%%%%%%%%%%%%%%%%%%%%%
The \SBD iterations are summarized in \eqref{eqs:SBD}.
\begin{equation}\label{eqs:SBD}
\left\{\begin{array}{l}
       \left\{\begin{array}{ll}
         \ m_i^{n+1} & =m^n_i  -\etaintwo m^n_i, \quad i\neq \imin:= \mathop{\argmin}_{i}F(\bxin) \\ \\
          \ m_\imin^{n+1} & =\displaystyle m_\imin^n +\sum \limits_{i\neq \imin}\etaintwo m_i^{n},  
     \end{array} \right\}  \etain := \displaystyle \Big(\frac{F(\bxin)-F^n_{\textnormal{min}}}{F^{n}_{\textnormal{max}}-F^{n}_{\textnormal{min}}}\Big)^\mq\\
     \\
       \displaystyle \qquad \relmin:= \frac{m^{n+1}_i}{m^{n+1}_+}, \quad m^{n+1}_+:=\max_i m^{n+1}_i \\ \\ 
     \left\{\begin{array}{rl}
     \bpin &:= \bpin(\bxin,\relmin) \  
    \left\{\begin{array}{l}
     \textnormal{Choose a random} \  r \in {\mathcal U} \big[\frac{1}{2}(1+\relmin), 1\big]; \\
     \textnormal{Algorithm} \ \ref{alg:random} \ \textnormal{computes}   \ \bpin \ \textnormal{such that}\\ 
     \left\langle \bpin,\nabla F(\bxin)\right\rangle =r|\nabla F(\bxin)|^2  
       \end{array}\right.
    \\ \\
    \hspace*{-0.2cm} \hin&:= h(\bxin,\clam\relmin) \quad \textnormal{Backtracking protocol in Algorithm} \ \ref{alg:backtracking}\\ \\
      \bx_i^{n+1} &= \bxin-\hin\bpin, 
     \end{array}\right.     
\end{array}\right. 
\end{equation}
The first part encodes the mass transfer from high to low ground in terms of 
 a communication protocol, that dictates mass transition factors, , $\{\etain\}$.
The second part encodes the stepping in a descent direction, $\hin\bpin$, based on two mass-dependent procedures: \newline
(i) a random choice of the descent direction, $\bpin=\bpin(\bxin,\relmin)$, whose  orientation is aligned within a random opening away from the orientation of the gradient $\nabla F(\bxin)/|\nabla F(\bxin)|$;  and\newline
(ii) a backtracking strategy for adjusting the step size, $\hin=h\big(\bxin,\clam\relmin\big)$, which secures the desired descent property. Observe that both the direction and step size are adjusted to the position and the relative mass of a given agent.\newline
These procedures are summarized in the following pseudo-codes.

 \subsection{A protocol for random choice of the descent direction}\label{sec:random}
Algorithm \ref{alg:random} picks a random orientation lying in the spherical cap of the unit sphere, $\boin\in \mathbb{S}^{d-1}$, centered around the gradient orientation, $\displaystyle \bqin=\frac{\nabla F(\bxin)}{|\nabla F(\bxin)|}$, and then sets the descent direction $\bpin=|\nabla F(\bxin)|\boin$.
To this end, we proceed in two steps. First,  sampling a randomly chosen point, $\bX=\big(X(1), \ldots, X(d-1), X(d)\big)\in \mathbb{S}^{d-1}$,  in the spherical cap centered around the north pole, $\bz=(0,0,\ldots, 1)$,
\[
 X(i)=\left\{\begin{array}{ll}
\displaystyle \sqrt{1-r^2}\frac{Y(i)}{|\bY|} & Y(i)\sim {\mathcal N}(0,1), \ i=1,2, \ldots d-1,\\
r, & i=d. \end{array}\right.
\]
Note that $\bY=\big(\sum_i Y^2(i)\big)^{-\nicefrac{1/}{2}}\big(Y(1),\ldots, Y(d-1)\big)$ is a random point (with normally distributed components) on $\mathbb{S}^{d-2}$ and therefore $\bX$ above is the projection of that random $\bY$  onto the spherical cap of ${\mathbb S}^{d-1}$ dictated by $r$. Here $r$ is a randomly chosen parameter from a uniform distribution in 
$ \frac{1}{2}(1+\relmin)< r < 1$; thus,  the spherical cap, shown as the shaded area in figure \ref{fig:spherical_cap}, has an opening angle of $\theta=\arccos(r)$,  ranging from $\theta=60^{\circ}$ for lightest agents to  the gradient orientation, $\theta=0^{\circ}$, for the heaviest agent.
 In the second step,   Algorithm \ref{alg:random} uses the unitary (Householder) reflection which reflects the north pole $\bz$ to    $\bqin$  
\[
{\mathbb P}_i^n={\mathbb I}-2\frac{\bvin(\bvin)^\top}{|\bvin|^2}, \qquad \bvin:=\bqin-\bz,
\]
and then reflects  $\bX$ into the desired  $\boin:= {\mathbb P}_i^n\bX$, see  Fig.~\ref{fig:spherical_cap},

\ifx%%
\begin{figure}[ht]
    \centering
    \includegraphics[width=0.4\textwidth]{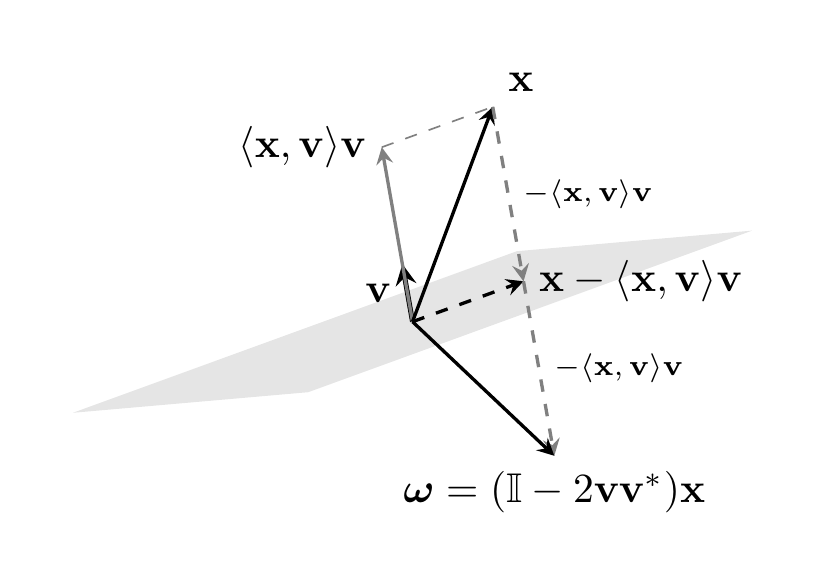}
    \caption{The classical (Householder) reflection efficiently rotates one vector towards the direction defined by another.}
\end{figure} 
\fi%%

\begin{algorithm}[htbp]
\begin{algorithmic}
\setstretch{1.1}
\State \caption{Random descent direction $\mathbf{p}_i^n$ for agent $\bxin$ with relative mass $\relmin$}
\label{alg:random}
\State Set $\displaystyle \mathbf{q}_i^n = \frac{\nabla F(\mathbf{x}_i^n)}{|\nabla F(\mathbf{x}_i^n)|}$
\State Choose random $r$ such that $\frac{1}{2}(1+\widetilde{m}_i^n) < r < 1$
\State Set random vector $\mathbf{Y} \in \mathbb{R}^{d-1}$ with $Y(i) \sim {\mathcal N}(0, 1)$ so that $\mathbf{Y}/|\mathbf{Y}| \in \mathbb{S}^{d-2}$
\For {$i=1$ to $d-1$}
\State Set $\displaystyle X(i) = \sqrt{1-r^2} \frac{Y(i)}{|\mathbf{Y}|}$
\EndFor
\State Set $X(d) = r$ so $\bX=(X(1), \cdots, X(d-1),X(d))\in {\mathbb S}^{d-1}$
\If {$1-\mathbf{q}_i^n(d) \neq 0$}
\State Set $\mathbf{v}_i^n = \mathbf{q}_i^n - \mathbf{z}$ where $\bz:=(0,\dots,0,1)$ is the north pole of $\mathbb{S}^{d-1}$
\State Set $\displaystyle \boin = \mathbf{X} - 2\frac{\langle \mathbf{v}_i^n, \mathbf{X} \rangle}{|\mathbf{v}_i^n|^2}\mathbf{v}_i^n$ \quad \ \  \% Simplification: {$|\bv_i^n|^2=2\big(1-\qin(d)\big)$}
\Else
\State Set $\boldsymbol{\omega}_i^n = \mathbf{X}$
\EndIf
\State Set $\mathbf{p}_i^n = |\nabla F(\mathbf{x}_i^n)| \boldsymbol{\omega}_i^n$
\end{algorithmic}
\end{algorithm}

  %%%%%%%%%%%%%%%%%%%%%%%%%%%%%%%%%%
 \subsection{Backtracking --- a protocol for time stepping}\label{sec:BT}
 %%%%%%%%%%%%%%%%%%%%%%%%%%%%%%%%%%%
The direction $\bpin$ computed in Algorithm \ref{alg:random} is partially aligned with $\nabla F(\bxin)$ so that \eqref{eq:porient} holds. Once the direction $\bpin$ is set, the new position $\bx_i^{n+1}(h)=\bxin-h\bpin$
is viewed as a function of
the step size $h$, and the objective is to  select an appropriate  step size, $\hin=h\big(\bxin,\clam\relmin\big)$, which  ensures the corresponding descent bound \eqref{eq:halfdescent},
\begin{equation}\label{eq:SBDdescent}
F(\bx^{n+1}_i) \leq F(\bxin)- \frac{1}{2}\clam\relmin \hin |\nabla F(\bxn)|^2, \qquad \quad \bx^{n+1}_i= \bxin-\hin\bpin, \quad 0<\lambda <1.
\end{equation}
A proper strategy for choosing such step size  is based on the classical
backtracking line search, \cite[\S3]{nocedal2006conjugate},  which is a computational realization of the well-known Wolfe
conditions \cite{wolfe1969,armijo1966}.
Recall that by Taylor's expansion \eqref{eq:halfdescent}, the desired  bound holds  for any  sufficiently small step size, $\hin\ll 1$. In fact, any step size $h <\nicefrac{1}{L}$, which need not be small, will suffice for the descent property, except that we do not have apriori access to the value of $L$. Our aim, therefore, is to choose a relatively large step size $h$, that even if  not optimally tuned with $\nicefrac{1}{L}$, it is still large enough to enforce the descent term $\hf \clam \relmin \hin|\nabla F(\bxin)|^2$.  To this end, one employs a
dynamic adjustment, starting with a relatively large $h=h_0$ (say -- $h_0=1$) for which one expects 
\[
F\big(\bxin-h\bpin\big)>F(\bxn)- \hf\clam \relmin h|\nabla F(\bxin)|^2,
\]
 and then successively shrink the step size,  $h\rightarrow \gamma h$, using a shrinkage factor $0<\gamma<1$, until the descent condition \eqref{eq:SBDdescent} is  fulfilled. Adjusting the shrinkage parameter $\gamma$ requires careful consideration of the trade-off between the  cost of a refined $\gamma \sim 1$ vs. improved performance with a crude $\gamma \ll 1$.\newline
The pseudo-code for computing the \SBD steps based on backtracking line search is given in Algorithm \ref{alg:backtracking} below.

\begin{algorithm}[ht!]
\setstretch{1.5}
\caption{Backtracking line search}
\label{alg:backtracking}
 \begin{flushleft}
Set the shrinkage parameter, $\gam\in(0, 1)$

\noindent
Set the relative mass $\displaystyle \relmin=\frac{m^{n+1}_i}{m^{n+1}_+}$

\noindent
Initialize the step size $h = h_0$.
\end{flushleft}

\begin{algorithmic}
\While{$F\big(\bxin-h\bpin\big)>F(\bxin)-\hf\clam\relmin h|\nabla F(\bxin)|^2$}

\State $h\gets \gamma h$.

\EndWhile

 \begin{flushleft}
Set $\hin=h\big(\bxin,\clam\relmin\big) \gets h$
\end{flushleft}
\end{algorithmic}

\end{algorithm}
\noindent

A stepping protocol for a non-convex optimization is required to strike a balance between small steps in  the vicinity of a potential minimizer and larger steps which avoid being trapped in local basins of attraction. The backtracking protocol achieves such a balance by adjusting the step size of each agent 
 according to its relative mass, $\relmin$ 
\begin{equation}\label{eq:step length}
    \hin=h\big(\bxin,\clam\relmin\big)  \qquad \relmin:=\frac{m^{n+1}_i}{m\mymax{n+1}}, \quad  m\mymax{n+1}:= \max_{i} m^{n+1}_i, \ \ 0<\clam<1,
\end{equation}
where $h(\bxin,\cdot)$ is a decreasing function of the relative mass $\clam\relmin$. 
Thus, our  mass-dependent backtracking is  an adaptive protocol:   it adapts itself from small time steps in the steepest gradient direction for heavier agents which lead the swarm, to larger steps in randomly chosen directions (that may differ from the steepest descent) for lighter agents which are the explorers of the swarm, exploring  the ambient space.

\medskip\noindent
{\bf The descent property}.  The backtracking \algo \ref{alg:backtracking} yields a step size $h_i^n=h(\bxin,\clam\relmin)$ with a lower bound
$\displaystyle  
\hin >  \frac{\gamma}{L}$ with $L$ denoting a Lipschitz bound of $\nabla F$ which is assumed to exists, $L=\max_{\bx\in \Omega} \|D^2F(\bx)\|<\infty$. Indeed, this can be argued by contradiction: if \mbox{$\displaystyle \frac{\hin}{\gamma} \leq  \frac{1}{L}$} then by \eqref{eq:porient} we would have, 
\[
\begin{split}
F\Big(\bxin-\frac{\hin}{\gamma}\bpin\Big)&\leq F(\bxn)-\frac{\hin}{\gamma}\langle \bpin,\nabla F(\bxin)\rangle+\frac{L}{2}\big(\frac{\hin}{\gamma}\big)^2|\bpin)|^{2}   \\
& \leq F(\bxin)-\Big(\frac{1+\widetilde{m}_i^n}{2}-\frac{L}{2}\frac{\hin}{\gamma}\Big)\frac{\hin}{\gamma}|\nabla F(\bxin)|^2\\
 &\leq F(\bxin)-\hf\clam\relmin \frac{\hin}{\gamma}|\nabla F(\bxin)|^2, \quad \clam<1.
\end{split}  
\]
But this contradicts the fact that  the backtracking iterations fail  to satisfy such inequality  with step size $\hin/\gamma$, since  according to \algo \ref{alg:backtracking}, $F(\bxin-(\hin/\gamma)\bpin)>F(\bxin)-\hf\clam\relmin (\hin/\gamma)|\nabla F(\bxin)|^2$ (in fact, the largest step size that succeeds in securing reverse inequality is with time step $\hin, \ \hin <\hin/\gamma$).
This contradiction confirms that $\ds \hin > \frac{\gamma}{L}$, which in turn enables us to convert the descent bound \eqref{eq:SBDdescent} into a precise descent property,
\begin{equation}\label{eq:SBDdescent-precise}
F(\bx_i^{n+1}) \leq F(\bxin)-\hf\clam\relmin \hin
|\nabla F(\bxin)|^2 \leq F(\bxin)-\frac{\gamma}{2L}\clam\relmin|\nabla F(\bxin)|^2
\end{equation}
The descent property we obtain is constrained by an additional factor of $\frac{1}{4}$ compared to the standard version of SBGD that relies on the gradient direction \cite[Proposition 5.2]{LTZ2022}. However, the randomization of the descent direction brings the advantage of allowing lighter agents to explore a wider range of directions. As we will see later, this exploration leads to substantial improvements in the optimization process in high dimensions.

It is important to note that heavy agents still adhere to the steepest descent along the gradient direction. The spherical cone of random directions is narrower for heavier agents. In fact, the heaviest agent strictly follows the steepest descent with $\bpin=\nabla F(\bxin)$, eliminating the need for a random choice at this particular point.

%%%%%%%%%%%%%%%%%%%%%%%%
\subsection{\SBD pseudocode}\label{sec:implement}
%%%%%%%%%%%%%%%%%%%%%%%%%
The pseudocode of the \SBD method is presented in Algorithm \ref{alg:SBD}. The initial setup involves $N$ randomly distributed agents $\bx^{0}_{1},\cdots,\bx_{N}^{0}$, associated with initial masses $m_{1}^{0},\cdots,m_{N}^{0}$. Initially, all  agents are assigned  equal masses, $\displaystyle m_{j}^{0} = \nicefrac{1}{N}$, $j = 1,\ldots, N$. At each iteration, the agent positioned at $\bx_\imin =  \argmin_{\bxin} F(\bxin)$ attains the minimal value, while the other agents transfer part of their masses to that minimizer $\bx_\imin$. Then all the agents are updated with the gradient descent method using the direction obtained in \eqref{eq:orient-omega} and step size in \eqref{eq:step length}. 

\begin{algorithm}[ht]
\setstretch{1.1}
\caption{Swarm-Based Random Descent Method}\label{alg:SBD}

\begin{algorithmic}
\State Set the parameters: $tolm$, $tolmerge$, $tolres$, and $nmax$
\State Set the number of agents,$N$, and the mass transfer parameter, $\mq\geq 1$
\State Randomly generate initial positions: $\bx_{1}^{0}, \ldots, \bx_{N}^{0}$
\State Set initial mass for all agents: $\displaystyle m^{0}_{1}=\cdots=m^{0}_{N} = \nicefrac{1}{N}$
\For{$n = 0, 1, 2, \ldots, nmax$}
\State Merge agents if their distance $< tolmerge$
\State Set the index of the optimal agent: $i_n = \argmin_{i} F(\bx^n_i)$
\State Set $F_{\textnormal{min}} = F(\bx^n_\imin)$ and $F_{\textnormal{max}} = \max_{i}F(\bxin)$
\For{$i = 1, \ldots, N$ \textbf{and} $i \neq \imin$} \hspace{0.6em} % Mass transitions
\If{$m_i^n < \nicefrac{1}{N} \cdot tolm$}
\State Set $m_i^{n+1} = 0$
\State Reduce the number of active agents: $N \gets N - 1$
\Else
\State Set $m_i^{n+1} = m^n_i - \etain m^n_i$ \quad where $\displaystyle \etain = \left(\frac{F(\bxin) - F^n_{\textnormal{min}}}{F^n_{\textnormal{max}} - F^n_{\textnormal{min}}}\right)^{\mq}$
\EndIf
\EndFor
\State Set $m_\imin^{n+1} = m_\imin^n + \sum_{i \neq \imin} \etain m^n_i$ \quad % The mass of the overall crowd is conserved
\State Set $m_+ = \max_i m^{n+1}_i$
\For{$i = 1, \ldots, N$} \hspace{0.7in} % Random-based descent iteration
\State Compute relative masses $\displaystyle \relmin=\frac{m^{n+1}_i}{m_+}$
\State Compute a random descent direction: $\bpin$ (using \algo \ref{alg:random})
\State Compute the step size: $h = h(\bxin, \clam\relmin)$ (using \algo \ref{alg:backtracking})
\State Update position: $\bx_i^{n+1} = \bxin - h \bpin$
\EndFor
\If{$|\bx_i^{n+1} - \bxin| \leq tolres$}
\State \textbf{break}
\EndIf
\EndFor
\end{algorithmic}
\end{algorithm}

\noindent
We use three tolerance factors:\newline
$\cdot$ $tolm$: If an agent's mass falls below this threshold, the agent is eliminated, and its remaining mass is transferred to the optimal agent at $\bx_\imin$.\newline 
$\cdot$ $tolmerge$: Agents that are sufficiently close to each other, i.e., their distance is below this threshold, are merged into a new agent. The masses of the merged agents are combined into the newly generated agent.\newline
$\cdot$ $tolres$: The iterations terminate when the descent of the minimizer between two consecutive iterations falls below this threshold. 

%%%%%%%%%%%%%%%%%%%
\section{Convergence and error analysis}\label{sec:convergence}
%%%%%%%%%%%%%%%%%%%%
The study of convergence and error estimates for the \SBD method requires quantifying the behavior of $F$. Here we emphasize that the required smoothness properties of $F$ are only sought in the region explored by the  \SBD iterations. We assume  that there exists a \emph{bounded} region, $\Omega \ni \bxin$ for all agents. Since the \SBD allows light agents to explore the ambient space with large step size (starting with $h_0$), we do not have an apriori bound on $\Omega$; in particular, the footprint of the \SBD crowd $\textnormal{conv}_i\{\bxin\}$ may expand well beyond its initial convex hull
$\textnormal{conv}_i\{\bx^0_i\}$. The expansion of the initial convex hull is an essential feature of the algorithm that allows the agents to find minima outside their initial range, demonstrated in the numerical experiments with shifted initial data domains such as in Table \ref{tab:SBD multiD off-centered ackley}. \newline
We consider the class of loss functions,  $F\in C^2(\Omega)$, with  Lipschitz bound $L=\max_\Omega\|D^2F\|<\infty$, 
\begin{equation}\label{eq:FisC2}
    |\nabla F(\bx)-\nabla F(\by)|\leq L|\bx-\by|,\quad \forall \bx,\by \in \Omega.
\end{equation}

\subsection{Convergence to a band of local minima} Our next proposition provides a precise quantitative description for the convergence of the \SBD method. The convergence   is determined by the time series of  \SBD minimizers, $\{\bXnmin\}$,
\begin{subequations}\label{eqs:XminXplus}
\begin{equation}\label{eq:Xmin}
\bXnmin=\bx^n_{i_n}, \qquad i_n:=\argmin_i F(\bxin).
\end{equation}
 We shall also need the time series of  its heaviest agents, $j_n:=\argmax_i m^n_i$; to this end, we let $\bX\mymax{n}$ denote the \emph{parent} of the heaviest agent at $t=t^{n+1}$
 \begin{equation}\label{eq:Xplus}
 \bX\mymax{n}=\bx^n_{j_{n+1}}, \qquad j_{n+1}:=\argmax_i m^{n+1}_i.
 \end{equation}
 \end{subequations}
The interplay between minimizers and the communication of masses leads to a gradual mass shift from higher ground to the minimizers. Eventually, the two sequences coincide when the \SBD minimizers gain enough mass to assume the role of heaviest agents.  
Finally, we  introduce the  scaling $M=\max_j F(\bx_j^0)-F(\bxmin)$ where $\bxmin$ is the global minimum. Since $F(\bxin)$ are decreasing, we conclude that the \SBD iterations remain within that range, namely
\begin{equation}\label{eq:var}
  F(\bx^n_i)-F(\bx^n_j) \leq M,  \qquad \forall n, j, \quad M:=\max F(\bx_i^0)-F(\bxmin).
 \end{equation}
\begin{proposition}\label{prop:convergence of SBGD}
Consider the \SBD iterations  \eqref{eqs:SBD} with  random-based search direction, $\bpin$, determined by \algo \ref{alg:random}, and with a step-size \eqref{eq:step length}, $\hin=h\big(\bxin,\clam\relmin\big)$, determined by backtracking line search of \algo \ref{alg:backtracking}.\newline
Let $\{\bXnmin\}_{n\geq0}$ and $\{\bX\mymax{n}\}_{n\geq0}$ denote the time sequence of \SBD  minimizers  and, respectively,   (parent of) heaviest agents outlined in \eqref{eqs:XminXplus}
 Then, there exists a constant, $C=C(\gamma,L,M, \clam)$ given in \eqref{eq:telesf} below, such that we have summability of gradients
\begin{equation}\label{eq:summability}
    \sum_{n = 0}^{\infty} \cdel_n^2\cdot \min\big\{1,  \cdel_n^{2\mq}\big\}  < C\min_i F(\bx^0_i), \qquad \cdel_n:=\min\{|\nabla F(\bX\mymax{n})|,|\nabla F(\bXnmin)|\}.
\end{equation}
Here, $\mq\geq 1$ is the mass transfer parameter in \eqref{eq:etai}. 
\end{proposition}

\noindent \emph{Proof.}  Our purpose is to find a lower bound on the relative masses, $\displaystyle \relmin=\frac{m^{n+1}_{i}}{m^{n+1}_{j_{n+1}}}$,
which will dictate the descent property of the different agents according to \eqref{eq:SBDdescent-precise}. 
Observe that for the heaviest agent, $i=j_{n+1}$, 
 \eqref{eq:SBDdescent-precise} with $\widetilde{m}^{n+1}_{j_{n+1}}=1$ implies
 \begin{equation}\label{eq:heaviest}
 F(\bx^{n+1}_{j_{n+1}}) \leq F(\bX^n_+)-\frac{\gamma}{2L}\clam|\nabla F(\bX^n_+)|^2, \qquad \bX^n_+=\bx^n_{j_{n+1}}.
 \end{equation}
We distinguish between two scenarios. The first is a canonical scenario in which the minimizing agent at $t=t^n$ coincides with the heaviest agent at time $t^{n+1}$, namely, when
$i_{n}=j_{n+1}$, or $\bXnmin=\bX^n_+$. Then \eqref{eq:heaviest}  implies
\begin{equation}\label{eq:telesa}
F(\bXnplusmin) \leq F(\bx^{n+1}_{j_{n+1}})  \leq  F(\bXnmin)-\frac{\gam}{2L}\clam|\nabla F(\bXnmin)|^2, \qquad \bXnmin=\bX^n_+.
\end{equation}
The inequality on the left follows since $\bXnplusmin=\bx^{n+1}_{i_{n+1}} $ is  the global minimizer at $t^{n+1}$.

\noindent
Next, we consider the second scenario  $i_{n}\neq j_{n+1}$, that is --- 
 when the mass of the minimizer $m^{n+1}_\imin$ did not yet `catch-up'  the position as the heaviest agent so that  $\displaystyle \widetilde{m}^{n+1}_{i_n}= \frac{m^{n+1}_\imin}{m^{n+1}_{j_{n+1}}}<1$. Yet, we claim that the descent property associated with the relative mass $\relminin$ cannot be arbitrarily small.
 We  consider two sub-cases, depending on the size of $\ds F(\bX\mymax{n})-F(\bXnmin)$.\newline
Case (i). Assume $\ds F(\bX\mymax{n})-F(\bXnmin) \leq \frac{\gam}{4L}\clam|\nabla F(\bX\mymax{n})|^2$.
Appealing to  \eqref{eq:heaviest}  we find
\begin{equation}\label{eq:telesc}
\begin{split}
F(\bXnplusmin) & \leq F(\bx^{n+1}_{j_{n+1}}) \leq
F(\bX\mymax{n}) -\frac{\gam}{2L}\clam|\nabla F(\bX\mymax{n})|^2 
  \leq F(\bXnmin) -\frac{\gam}{4L}\clam|\nabla F(\bX\mymax{n})|^2.
\end{split}
\end{equation}
The inequality on the left follows since $\bXnplusmin$ is  the global minimizer at $t^{n+1}$; the middle inequality quotes \eqref{eq:heaviest}
and the last inequality follows from our assumption.\newline
Case (ii). Finally, we remain with the case
\[
F(\bX\mymax{n})-F(\bXnmin) > \frac{\gam}{4L}\clam|\nabla F(\bX\mymax{n})|^2.
\]
We claim that in this case,  
\begin{equation}\label{eq:claimedm}
\relminin >\frac{1}{M^2}\big(F(\bX\mymax{n})-F(\bXnmin)\big)^2 \geq  \Big(\frac{\gam\clam}{4ML}\Big)^{\mq}|\nabla F(\bX\mymax{n})|^{2\mq}.
\end{equation}
Indeed,  since agent $j_{n+1}$ is not the minimizer at time $t=t^n$, namely $ j_{n+1}\neq i_n$,  
then it had to shed a  portion of its mass, $ m^n_{j_{n+1}}-\eta\mymax{n} m^n_{j_{n+1}} \rightarrow m^{n+1}_{j_{n+1}}$, which was transferred to the 
minimizer $m^{n+1}_\imin \leftarrow m^n_\imin + \ldots + \eta\mymax{n}m^n_{j_{n+1}}$. Thus, the loss of mass by heavy agent 
\[
m^{n+1}_{j_{n+1}}=m^n_{j_{n+1}}-\eta\mymax{n} m^n_{j_{n+1}},\quad
\eta\mymax{n} =  \Big(\frac{F(\bx^n_{j_{n+1}})-F(\bx^n_{i_n})}{\max_j F(\bx^n_j)- F(\bx^n_\imin)}\Big)^{\mq} \geq \frac{1}{M^{\mq}} \big(F(\bX^n_+)-F(\bXnmin)\big)^{\mq}.
\]
was \emph{gained}  by the minimizer agent, $i=\imin$. Therefore, the relative mass of that minimizer is at least as large as claimed in \eqref{eq:claimedm}
\[
 \relminin = \frac{m^{n+1}_\imin}{m^{n+1}_{j_{n+1}}}> \frac{\eta\mymax{n}}{1-\eta\mymax{n}} \geq \frac{1}{M^{\mq}} \big(F(\bX^n_+)-F(\bXnmin)\big)^{\mq} \geq \Big(\frac{\gam\clam}{4ML}\Big)^{\mq}|\nabla F(\bX\mymax{n})|^{2\mq}.
\] 
The descent property \eqref{eq:SBDdescent-precise} together with \eqref{eq:claimedm}  imply
\begin{equation}\label{eq:telesb}
\begin{split}
F(\bXnplusmin)    & \leq F(\bx^{n+1}_{i_n}) \leq F(\bx^n_{i_n}) - \frac{\gam}{2L}\clam\relminin|\nabla F(\bx^n_{i_n})|^{2}
 \\
& \leq F(\bXnmin)-\frac{\gam\clam}{2L}\Big(\frac{\gam\clam}{4ML}\Big)^{\mq}|\nabla F(\bX\mymax{n})|^{2\mq}\cdot|\nabla F(\bXnmin)|^2.
\end{split}
\end{equation}
Combining \eqref{eq:telesa}, \eqref{eq:telesc} and \eqref{eq:telesb}  we find 
\begin{equation}\label{eq:telesf}
\begin{split}
F(\bXnplusmin) & \leq F(\bXnmin) -\frac{1}{C} \min\big\{|\nabla F(\bXnmin)|^{2}, |\nabla F(\bX\mymax{n})|^{2}, |\nabla F(\bX\mymax{n})|^{2\mq}\cdot|\nabla F(\bXnmin)|^{2}\big\}\\
& \leq F(\bXnmin) -\frac{1}{C} \cdel_n^2 \min\big\{1, \cdel_n^{2\mq}\big\}, \qquad  C=\frac{4L}{\gamma\clam}\cdot\max\Big\{2,\Big(\frac{4ML}{\gamma\clam}\Big)^{\mq}\Big\}.
\end{split}
\end{equation}
The desired bound \eqref{eq:summability} follows by a telescoping sum. $\square$

 The  summability bound \eqref{eq:summability} implies that eventually, for large enough $n>N_0$, the minimizers and (parent of) heaviest \SBD agents, $\cdel_n<1$ and hence
\[
\sum_{n>N_0}^\infty  \min\{|\nabla F(\bX\mymax{n})|,|\nabla F(\bXnmin)|\}^{2(\mq+1)} \leq C\min_i F(\bx^0_i).
 \] 
 It follows that there exist sub-sequences, $\bX^{n_\alpha}\in \{\bX\mymax{n}\}_{n\geq N_0}\cup \{\bXnmin\}_{n\geq N_0}$, satisfying 
the  Palais-Smale condition, \cite{palais1964generalized}, $F(\bX^{n_\alpha})\leq \max_i F(\bx^0_i)$ while $\nabla F(\bX^{n_\alpha}) \stackrel{\alpha\rightarrow \infty}{\longrightarrow}0$. Arguing along \cite[Theorem 5.4]{LTZ2022} we summarize by stating the following.

\begin{theorem}\label{thm:SBD convergence}
Let $\{\bXn\}_{n\geq N_0}:=\{\bX\mymax{n}\}_{n\geq N_0}\cup \{\bXnmin\}_{n\geq N_0}$ denote the combined time sequence of \SBD  minimizers/heaviest agents, \eqref{eqs:XminXplus}.   Then there exist   one or more sub-sequences, $\{\bXna, \ \alpha=1,2,\ldots,\}$, that converge to a band of local minima with equal heights,
\begin{equation}\label{eq:convergence}
\bXna \stackrel{\alpha\rightarrow \infty}{\longrightarrow}\bX^*_\alpha \ \ \textnormal{such that} \  \nabla F(\bX^*_\alpha)=0, \ \textnormal{and} \ F(\bX^*_\alpha)=F(\bX^*_\beta)
\end{equation}
In particular,in the generic case that $F$ admits only distinct local minima in $\Omega$, namely --- different local minima have different heights, then the \underline{whole} sequence $\bX^n$ converges to a local minimum.
\end{theorem}
\noindent \emph{Proof.} Since we assume  the sequence $\{\bXnmin\}$ is bounded in $\Omega$, it has  a converging sub-sequences. Take \emph{any} such converging sub-sequence $\bXnmina \rightarrow \bX^*_\alpha \in \Omega$. By \eqref{eq:summability}, $\nabla F(\bXnmina) \rightarrow 0$ for all sub-sequences, and hence $\bX^*_\alpha$ are local minimizers, $\nabla F(\bX^*_\alpha)=0$. Moreover, since $F(\bXnmin)$ is a decreasing, 
all $F(\bX^*_\alpha)$ must have the same `height'. The collection of equi-height minimizers  $\{\bX^*_\alpha \ \big| \  F(\bX^*_\alpha)=F(\bX^*_\beta)\}$ is the limit-set of $\{\bXnmin\}$. $\square$

Moreover, for analytic $F$'s, we can quantify the convergence \emph{rate} \eqref{eq:convergence}. To this end we use   Lojasiewicz inequality, \cite{law1965ensembles,lojasiewicz1993geometrie}, which guarantees that each critical point of analytic $F$ has ``flatness'' of some fixed order $\beta\in (1,2]$ in the sense  that there exists a neighborhood $\calN_* \ni\bxmin$ surrounding $\bxmin$, an exponent  $\beta$ and a constant $\mu>0$ such that
\begin{equation}\label{eq:Lojasiewicz}
\mu|F(\bx)-F(\bxmin)| \leq |\nabla F(\bx)|^\beta, \qquad \forall \bx\in \calN_*.
\end{equation}

\begin{theorem}\label{thm:GD-backtracking with PL}
 Consider an analytic loss function $F$ with minimal flatness 
 $\beta \in(1,2]$, such that the Lipschitz bound \eqref{eq:FisC2} holds.  Let $\{\bXnmin\}_{n\geq0}$ denote the time sequence of   \SBD minimizers,  \eqref{eqs:SBD},\eqref{eq:step length}, with converging sub-sequence, $\{\bXnmina\}_\alpha$, outlined in Theorem \ref{thm:SBD convergence}.  Then, there exists a constant, $C=C(\gamma,\lambda,\mu)$, such that 
\begin{equation}\label{eq:GD PL convergence}
    F(\bXnmina)-F(\bX^*_\alpha) \lesssim 
   \Big( \frac{C}{n_\alpha}\Big)^{\beta'}, \qquad \beta'=\frac{\beta}{2(\mq+1)-\beta}, \ \  \beta\in (1,2).
 \end{equation}
\end{theorem}

\noindent
Observe that as `flatness', increases, $\beta$ decreases the polynomial decay in \eqref{eq:GD PL convergence}. A more careful analysis which we omit\footnote{Requires to eliminate case (ii) in the proof of proposition \ref{prop:convergence of SBGD}; consult \cite{LTZ2022}},  allows to replace the factor $(\mq+1)$ by $\mq$, in which case, \eqref{eq:GD PL convergence} with $\mq=1, \beta=2$ implies exponential convergence. 

\medskip\noindent 
\emph{Proof.} We summarize the different statements  of  descent properties in  \eqref{eq:telesa}, \eqref{eq:telesc} and \eqref{eq:telesb}, writing
\[
F(\bX^{n+1}_-)\leq 
F(\bX^n_-) -\frac{1}{C}|\nabla F(\bX^n_\pm)|^{2(\mq+1)}, \qquad n>N_0,
\]
where $\bX^n_\pm=\argmin_{\bX^n_\pm} \{|\nabla F(\bX^n_+)|,|\nabla F(\bX^n_-)|\}$.
We focus on the converging sub-sequence $\{\bXnmina\}$,
\begin{equation}\label{eq:subseq}
F(\bX^{n_\alpha+1}_-)-F(\bX^*_\alpha)\leq F(\bX^{n_\alpha}_-)-F(\bX^*_\alpha) -\frac{1}{C}|\nabla F(\bX^{n_\alpha}_\pm)|^{2(\mq+1)}.
\end{equation}
Using Lojasiewicz bound \eqref{eq:Lojasiewicz},  and the fact that $F(\bX^n_+)\geq F(\bX^n_-)$, we find
\begin{equation}\label{eq:Loj}
|\nabla F(\bX^{n_\alpha}_\pm)|^\beta \geq \mu|F(\bX^{n_\alpha}_\pm)-F(\bX^*_\alpha)| \geq \mu|F(\bX^{n_\alpha}_-)-F(\bX^*_\alpha)|.
\end{equation}
Combining \eqref{eq:subseq}, \eqref{eq:Loj}, we conclude that the   error, $E_{n_\alpha}:=F(\bX^{n_\alpha}_-)-F(\bX^*_\alpha)$, satisfies
 \[
   E_{n_\alpha+1} \leq E_{n_\alpha}-\frac{1}{C}
   (\mu E_{n_\alpha})^{\frac{2(\mq+1)}{\beta}},\qquad 
\bXnmina \in \calN_\alpha.
\]
The solution of this Riccati inequality  yields
\[
F(\bXnmina)-F(\bX^*_\alpha)\leq \left\{|\min_i F(\bx^0_i)-F(\bX^*_\alpha)|^{-\nicefrac{1}{\beta'}}+ \frac{1}{C}\mu^{\frac{2(\mq+1)}{\beta}} n_\alpha\right\}^{-\beta'}, \quad \beta'=\frac{\beta}{2(\mq+1)-\beta}.
\]
and \eqref{eq:GD PL convergence} follows. $\square$

%%%%%%%%%%%%%%%%%%%%%
\section{Numerical results}\label{sec:results}
%%%%%%%%%%%%%%%%%%%%%%%

 Initially, the agents are placed  at random positions, $\{\bx^0_i\}$ with  equi-distributed masses $\{m^0_i=\nicefrac{1}{N}\}$. Masses are transferred from the high to the lowest ground at each iteration. Since we implement a ``survival of the fittest'' protocol in which the agent with the worst (=highest) configuration is eliminated, the swarm size decreases, one agent at a time, until only the heaviest agent remains. Our choice for the time-stepping protocol, $h\big(\bx,\clam\widetilde{m}\big)$, is the \emph{backtracking line search} outlined in \S\ref{sec:BT}, which is  weighted by the relative masses, $\relmin$.
The backtracking enforces a descent property for the \SBD iterations $\bxin$,
and the parameter,  $\clam\in (0,1)$, dictates how much the descent property holds in the sense that \eqref{eq:descentclam} is fulfilled.

\medskip\noindent 
We illustrate the performance of the multi-dimensional \SBD algorithm, \eqref{eqs:SBD},\eqref{eqs:orient-omega}, in several benchmark test cases \cite{benchmarks}.
The results are based on $k=1000$ runs of uniformly generated initial data in a hypercube. Backtracking parameters in Algorithm \ref{alg:backtracking} are $\clam=0.2$ and $\gamma=0.9$ and $h_0=1.$ The parameters in Algorithm \ref{alg:SBD} are $tolm = 10^{-4}$, $tolmerge = 10^{-3}$, $tolmax = 10^{-4}$ and $nmax=200$.\newline
We use the \emph{success rate} among the $k$ independent simulations to evaluate the solution's quality. We consider a simulation to be successful if $ \bx_{SOL}$ is within the $d$-dimensional ball of the global minimum: $ |\bxmin - \bx_{SOL} | \leq 0.1$. This condition ensures that the approximate solution lies in the basin of attraction of the global minimizer. 
In section \ref{sec:mqeq2} we fix the mass transfer parameter $\mq=2$; the effect of increasing $\mq=4,8$ is discussed in section \ref{sec:high-q}.

\subsection{Examples of \SBD with mass transfer parameter $\mq=2$}\label{sec:mqeq2} 

\begin{figure}[h]
    \centering
    \includegraphics[width=0.45\textwidth]{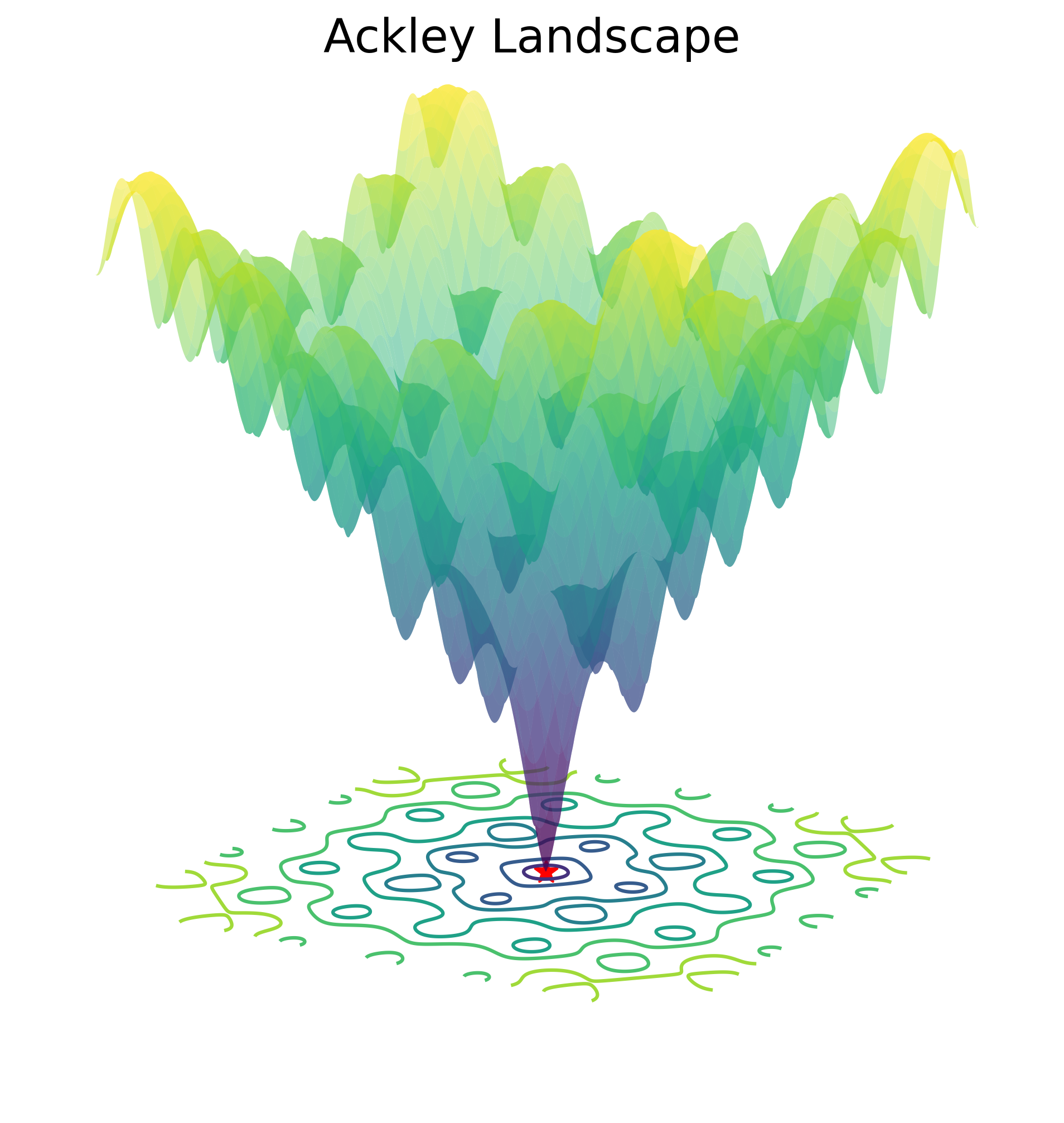}
    \includegraphics[width=0.45\textwidth]{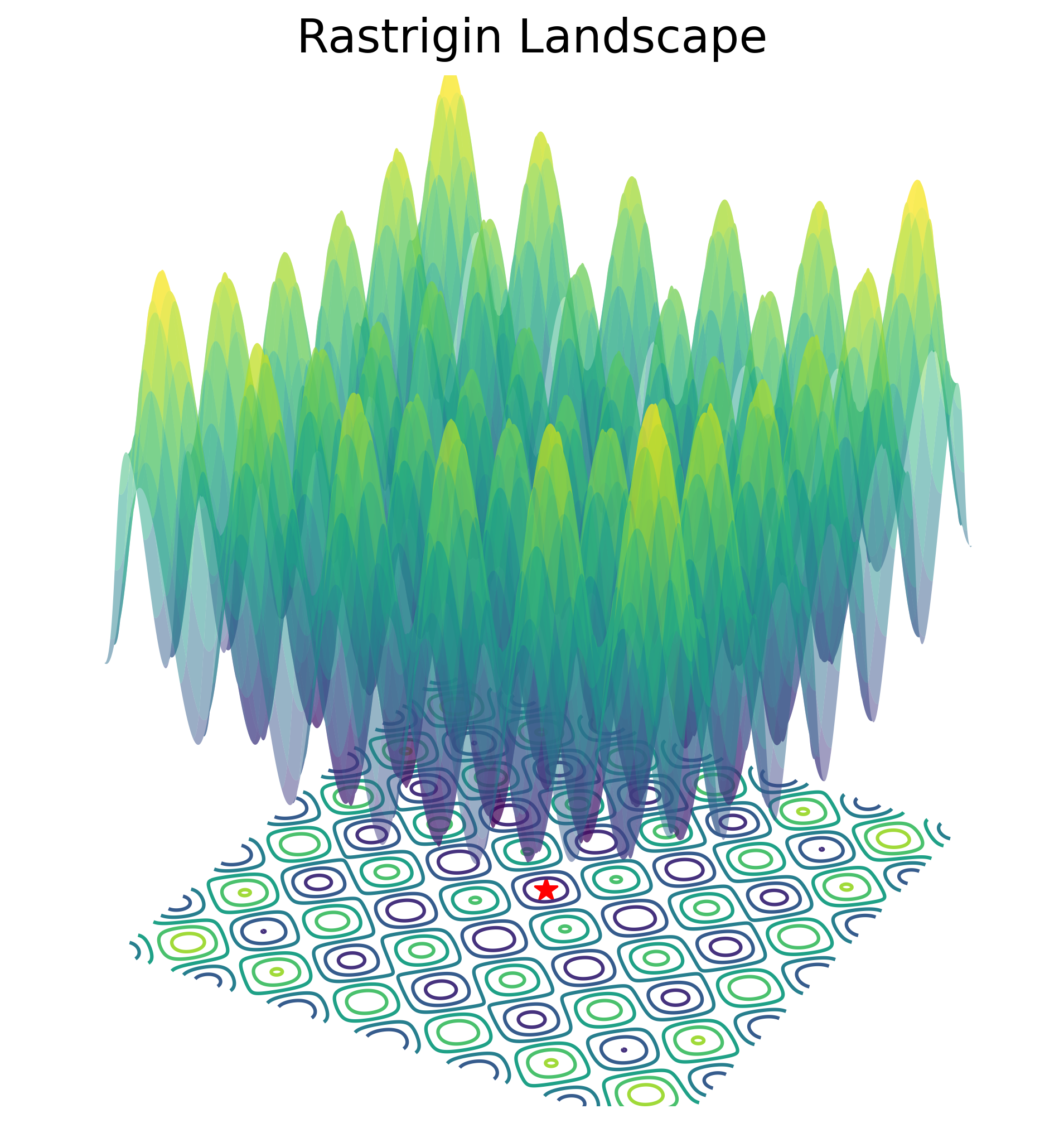}
    \includegraphics[width=0.45\textwidth]{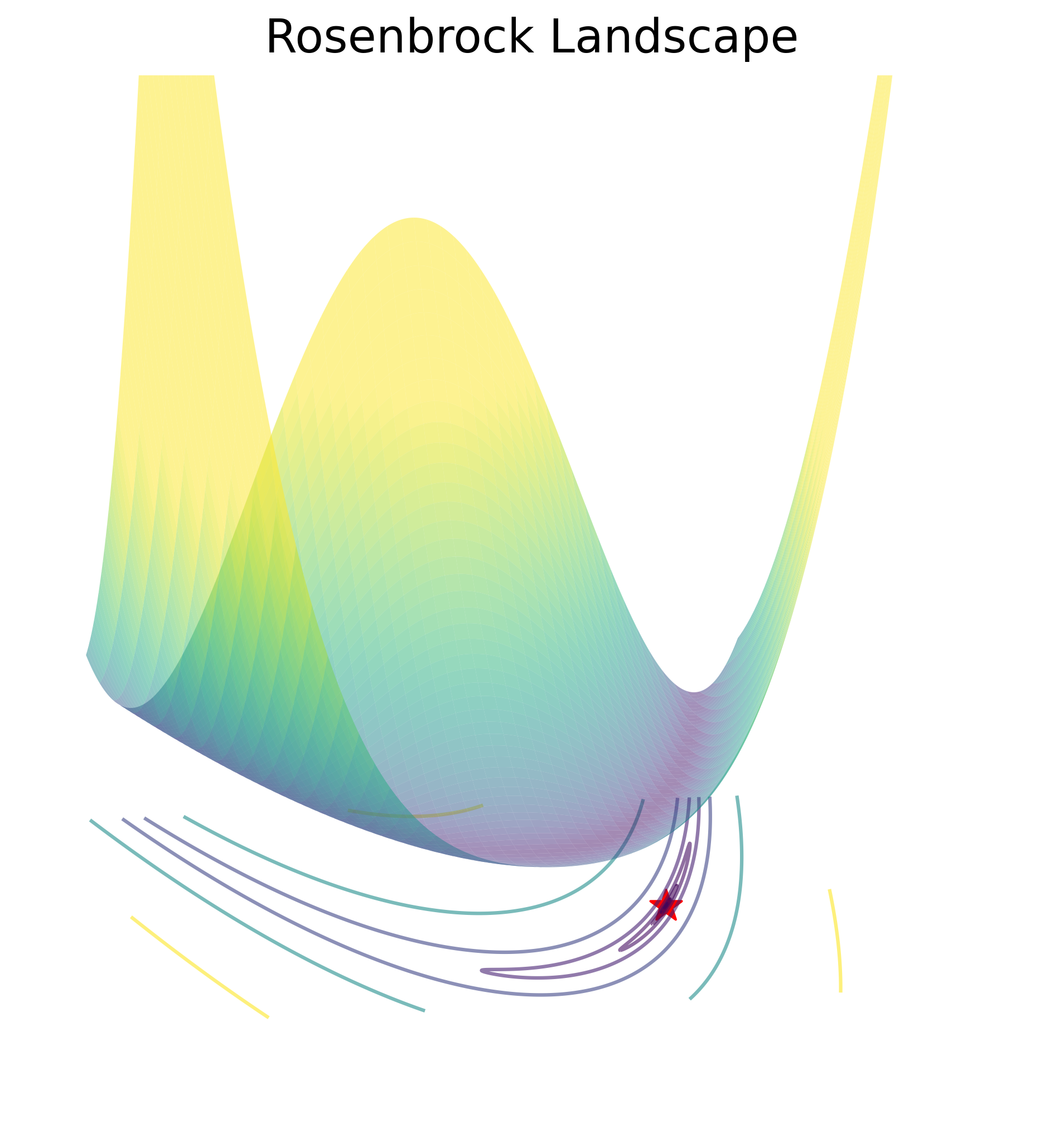}
    \includegraphics[width=0.45\textwidth]{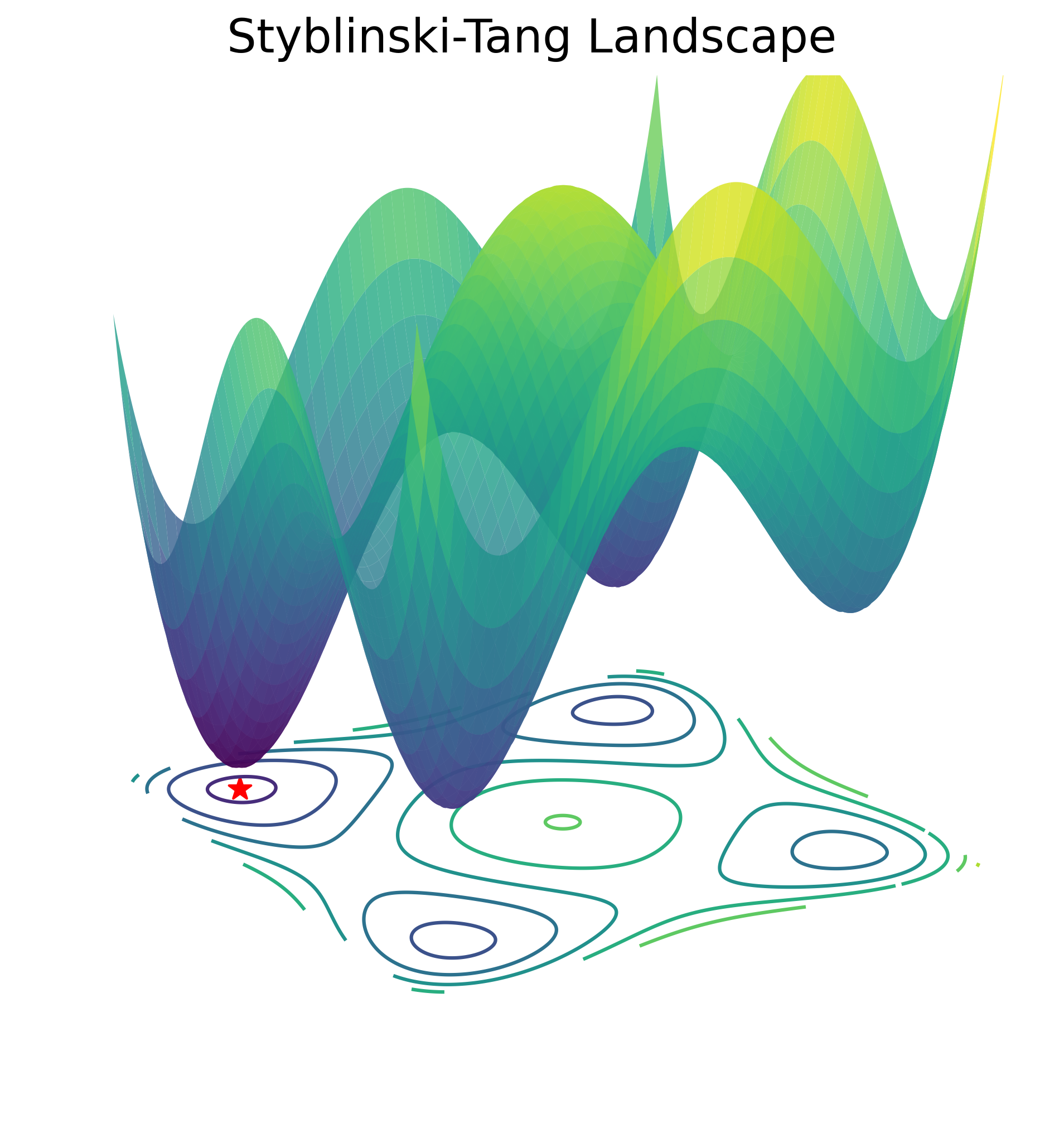}
       \caption{Two-dimensional landscapes for the test functions Ackley \eqref{eq:Ackley}, Rastrigin \eqref{eq:Rastrigin}, Rosenbrock \eqref{eq:Rosenbrock}, and Styblinski-Tang \eqref{eq:ST} with a contour plot on the bottom and a red star indicating the global minimum.}
       \label{fig:landscapes}
\end{figure}

Extensive comparisons of the gradient-based deterministic SBGD were performed in \cite{LTZ2022}. In this paper, we focus on the impact of randomization on success rates in comparison to SBGD. We consider four benchmarks using the Ackley, Rastrigin, Rosenbrock, Styblinski-Tang objective functions in $d$-dimensions.

\noindent
The \textbf{Ackley} function
\begin{equation}\label{eq:Ackley}
    F_{\textnormal{Ackley}}(\bx) = -20\exp\Big\{-\frac{0.2}{\sqrt{d}}\Big\{\sum^{d}_{i=1} x_i^{2}\Big\}^{\nicefrac{1}{2}}\Big\}-\exp\Big\{\frac{1}{d}\sum^{d}_{i=1}\cos(2\pi x_i)\Big\}+20+e,
\end{equation}
and the \textbf{Rastrigin} function
\begin{equation}\label{eq:Rastrigin}
    F_{\textnormal{Rstgin}}(\bx) = 10 d + \sum^{d}_{i=1}\Big\{ x_i^{2}-10\cos(2\pi x_i)\Big\}.
\end{equation}
have their global minimum at the origin, $\bxmin = 0$. The \textbf{Rosenbrock} function 
\begin{equation}\label{eq:Rosenbrock}
F_{\textnormal{Rsnbrk}}(\bx) = \sum_{i=1}^{d-1}\big(100(x_{i+1}-x_i^2)^2+(1 - x_i)^2\big),
\end{equation}
has its global minimum at $\mathbf{x}^* = (1, \ldots, 1)$. And finally, the \textbf{Styblinski-Tang} function
\begin{equation}\label{eq:ST}
F_{\textnormal{ST}}(\bx)=\frac{1}{2}\sum_{i=1}^d(x_i^4-16x_i^2 +5x_i),
\end{equation}
has its global minimum at $\mathbf{x}^* = (-2.903534, \ldots, -2.903534)$.
The two-dimensional landscapes of these benchmark examples are shown in Figure \ref{fig:landscapes}.

\begin{table}[h]
\setstretch{1.4}
\centering
{%\footnotesize
\begin{tabular}{|c| |cc|cc|cc|cc|cc|cc|}
\hline
\backslashbox{$d$}{$N$} & \multicolumn{2}{c|}{10} & \multicolumn{2}{c|}{25} & \multicolumn{2}{c|}{50} & \multicolumn{2}{c|}{100} \\ %& \multicolumn{2}{c|}{200} \\
\hline
& \SBD & SBGD & \SBD & SBGD & \SBD & SBGD & \SBD & SBGD \\ %& \SBD & SBGD\\
\hline
2 & \textbf{31.9\%} & 28.0\% 
  & \textbf{96.8\%} & 67.8\% 
  & \textbf{100.0\%}& 95.3\% 
  & 100.0\%& 100.0\%\% \\
%  & 100.0\%& 100.0\%\%\\

3 & 5.2\%& 5.6\% 
  & \textbf{17.6\%}& 13.6\% 
  & \textbf{57.9\%}& 28.6\% 
  & \textbf{92.4\%}& 52.0\% \\
%  & \textbf{99.8\%}& 77.3\%\\

4 & 0.3\% & 1.0\% 
  & 2.2\% & 3.9\% 
  & \textbf{7.2\%}  & 5.7\% 
  & \textbf{17.9\%} & 11.4\% \\
%   & \textbf{39.3\%} & 18.0\%\\

5 & 0.1\% & 0.0\% 
  & 0.2\% & 0.4\% 
  & 0.8\% & 0.4\% 
  & \textbf{3.2\%} & 1.2\% \\
%   & \textbf{4.1\%} & 2.9\%\\

6 & 0.0\% & 0.0\% 
  & 0.0\% & 0.0\% 
  & 0.0\% & 0.1\% 
  & 0.2\%& 0.4\%\\
%   & 0.6\%& 0.8\%\\
\hline
\end{tabular}
}
\smallskip
\caption{Success rates of \SBD vs. SBGD for global optimization of the $d$-dimensional Rastrigin function \eqref{eq:Rastrigin}.}
\label{tab:rastrigin}
\end{table}

\begin{table}[h]
\setstretch{1.5}
       \centering
    \begin{tabular}{|c| |cc|cc|cc|cc|cc|}
    \hline
    \backslashbox{$d$}{$N$} & \multicolumn{2}{c|}{10} & \multicolumn{2}{c|}{25} & \multicolumn{2}{c|}{50} & \multicolumn{2}{c|}{100} \\
\hline
& \SBD & SBGD & \SBD & SBGD & \SBD & SBGD & \SBD & SBGD\\
\hline
2 & \textbf{12.0\%} & 10.3\% & \textbf{52.1\%} & 18.7\% & \textbf{92.7\%} & 39.4\% & \textbf{99.2\%} & 56.7\% \\
3 & 2.4\% & 2.2\% & 8.1\% & 9.6\% & 27.2\% & 33.9\% & \textbf{82.6\%} & 71.0\% \\
4 & 2.3\% & 2.1\% & 3.5\% & 3.0\% & \textbf{9.4\%} & 3.9\% & \textbf{27.0\%} & 6.5\% \\
5 & 1.1\% & 0.8\% & 1.3\% & 1.6\% & \textbf{5.9\%} & 3.2\% & \textbf{10.2\%} & 6.1\% \\
6 & 0.5\% & 0.6\% & 1.1\% & 1.2\% & 1.6\% & 1.7\% & \textbf{5.1\%} & 2.6\% \\
\hline
\end{tabular}
\smallskip
\caption{Success rates for the optimization of the Rosenbrock function \eqref{eq:Rosenbrock} with initial agents $\bx_i^0\in[-2.048,2.048]^d$.}
\label{tab:rosenbrock}
\end{table}

\begin{table}[h]
\vspace*{-0.1cm}
\setstretch{1.5}
       \centering
    \begin{tabular}{|c| |cc|cc|cc|cc|cc|}
    \hline
    \backslashbox{$d$}{$N$} & \multicolumn{2}{c|}{10} & \multicolumn{2}{c|}{25} & \multicolumn{2}{c|}{50} & \multicolumn{2}{c|}{100} \\
\hline
& \SBD & SBGD & \SBD & SBGD & \SBD & SBGD & \SBD & SBGD\\
\hline
2 & \textbf{97.0\%} & 92.8\% & 100.0\% & 99.9\% & 100.0\% & 100.0\% & 100.0\% & 100.0\% \\
%3 & 55.3\% & 62.9\% & \textbf{99.1\%} & 96.3\% & 100.0\% & 99.7\% & 100.0\% & 100.0\% \\
4 & 29.5\% & 35.3\% & \textbf{83.7\%} & 79.0\% & \textbf{99.2\%} & 97.4\% & 100.0\% & 99.9\% \\
%5 & 16.4\% & 19.4\% & 47.4\% & 55.0\% & \textbf{86.0\%} & 83.6\% & \textbf{99.0\%} & 97.3\% \\
6 & 7.8\% & 10.4\% & 28.5\% & 32.5\% & 54.5\% & 55.4\% & \textbf{86.3\%} & 83.2\% \\
%\hline
%7 & 3.7\% & 5.2\% & 15.7\% & 19.2\% & 32.9\% & 34.2\% & \textbf{60.6\%} & 58.5\% \\
%\hline
8 & 2.2\% & 2.5\% & 7.6\% & 9.7\% & 13.7\% & 18.7\% & \textbf{36.7\%} & 35.4\% \\
%\hline
%9 & 1.6\% & 2.2\% & 3.7\% & 4.9\% & 8.5\% & 11.0\% & 20.0\% & 21.4\% \\
%\hline
10 & 0.4\% & 0.6\% & 2.6\% & 3.2\% & 5.9\% & 6.0\% & 10.2\% & 12.5\% \\
12 & 0.1\% & 0.2\% & 0.5\% & 0.8\% & 1.3\% & 2.2\% & 2.9\% & 3.8\% \\
\hline
\end{tabular}
\smallskip
 \caption{Success rates for the optimization of the Styblinski-Tang function \eqref{eq:ST} with initial agents $\bx_i^0\in[-3,3]^d$.}
 \label{tab:styblinski-tang}
 \end{table}

\medskip
Tables \ref{tab:rastrigin}, \ref{tab:rosenbrock} and \ref{tab:styblinski-tang} show the  advantage of \SBD over SBGD for  Rastrigin,  Rosenbrock and Styblinski-Tang functions.
We bold-face the success rate of \SBD in the tables when the advantage of \SBD over SBGD is at least 1\%. Observe that the advantage of randomization is only relevant in higher dimensions where SBGD has a very low success rate.
We recall that the same improved success rate of \SBD over SBGD was already recorded for the Ackley test function in Table \ref{tab:SBD multiD ackley}.
 This observation remains  valid when the range of initial agents for Ackley test function lies outside the neighborhood of its global minimum; this is documented in Table \ref{tab:SBD multiD off-centered ackley}. 

\begin{table}[h]
%\vspace*{-0.7cm}
\setstretch{1.5}
       \centering
    \begin{tabular}{|c| |cc|cc|cc|cc|cc|}
    \hline
    \backslashbox{$d$}{$N$} & \multicolumn{2}{c|}{10} & \multicolumn{2}{c|}{25} & \multicolumn{2}{c|}{50} & \multicolumn{2}{c|}{100} \\
\hline
& \SBD & SBGD & \SBD & SBGD & \SBD & SBGD & \SBD & SBGD\\
\hline
12 & 2.8\% & 3.7\% & 39.3\%&60.1\% & 74.8\%&96.2\% & 94.5\%&99.9\%\\
14 &  0.3\%&0.0\% & \textbf{19.6\%} & 0.9\% & \textbf{51.3\%}&2.0\% & \textbf{81.3\%} & 9.9\%\\
16 &  0.0\%&0.0\% & \textbf{2.7\%} & 0.0\% & \textbf{21.9\%}&0.0\% & \textbf{47.4\%} & 0.0\%\\
18 &  0.0\%&0.0\% & 0.0\% & 0.0\% & 0.7\%&0.0\% & \textbf{7.3\%} & 0.0\%\\
       \hline
        \end{tabular}
         \smallskip
 \caption{Same as Table \ref{tab:SBD multiD ackley} except the range of initial agents with $\bx_i^0\in[-3,-1]^d$ does not contain the global minimum of the Ackley function.}
 \label{tab:SBD multiD off-centered ackley}
 \end{table}

\begin{figure}[h]
\vspace*{-0.2cm}
    \centering
    \includegraphics[width=0.9\textwidth]{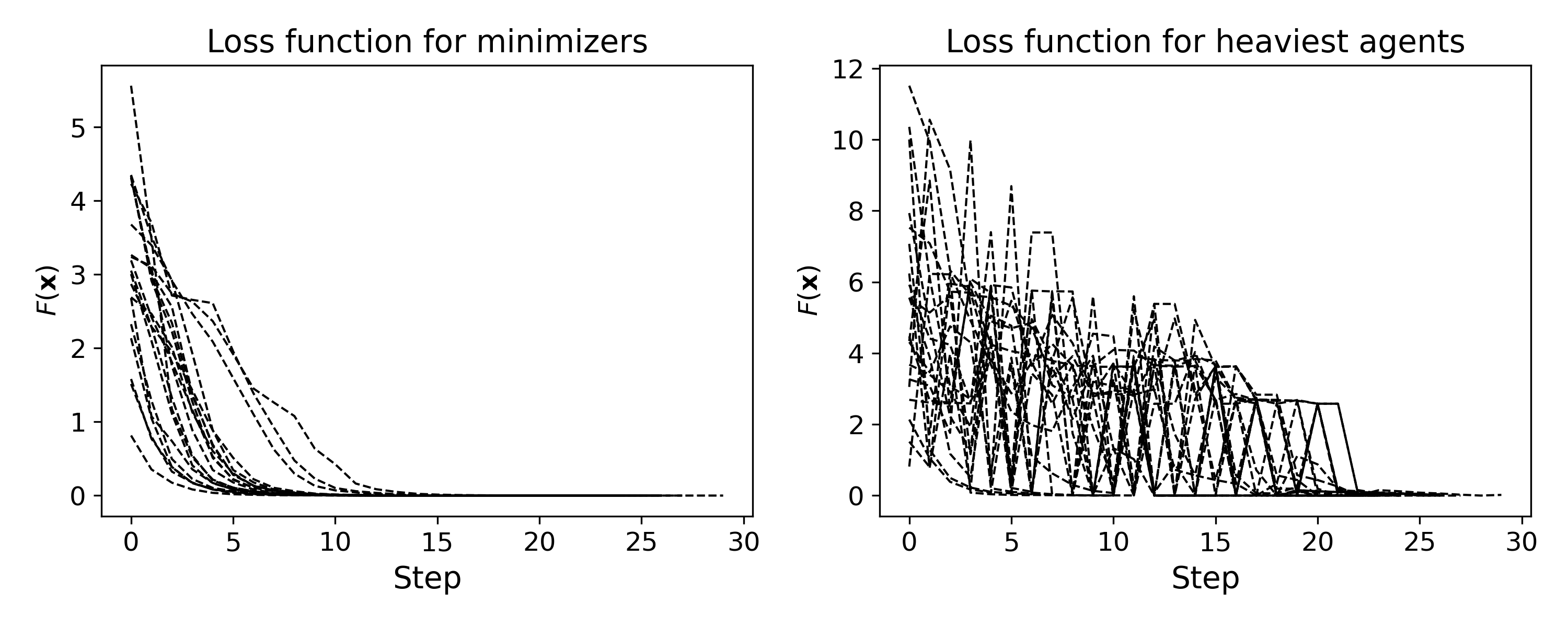}
       \caption{Loss functions for minimizers and heaviest agents as defined in \eqref{eqs:XminXplus} and \eqref{eq:Xplus} during optimization of the two-dimensional Ackley function for $20$ simulations with $N=50$ agents. Note the different scales on the $y$-axis between the minimizers and heaviest agents.}
       \label{fig:minimizers}
\end{figure}

\begin{figure}[h]
%\vspace*{-0.7cm}
    \centering
    \includegraphics[width=0.45\textwidth]{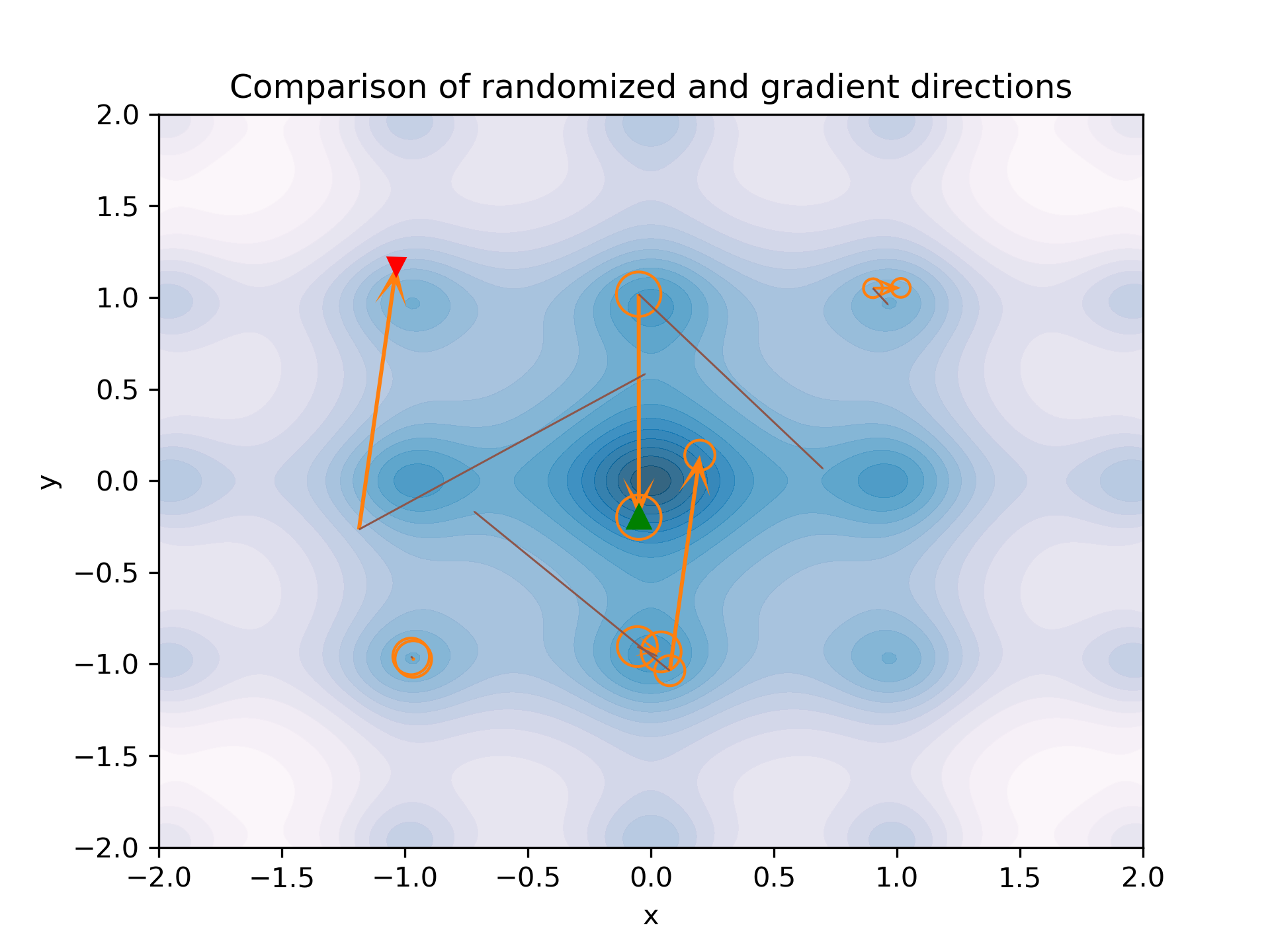}
    \includegraphics[width=0.45\textwidth]{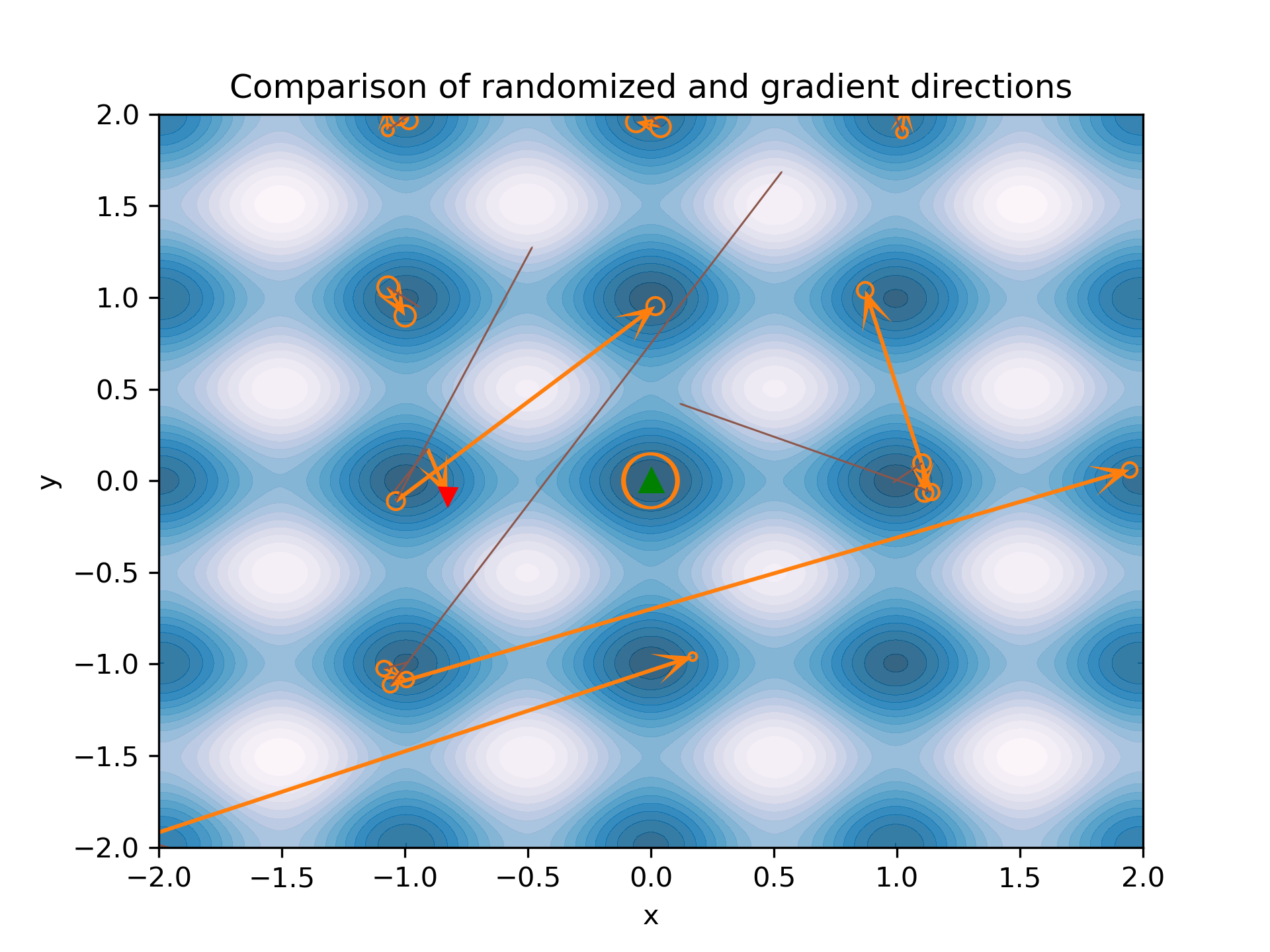}       
       \caption{A comparison of random descent, $-\bpin$ (orange arrows), and gradient descent, $-\nabla F(\bxin)$ (brown lines), during a simulation for the optimization of Ackley function (left) and the Rastrigin function (right).  The green triangle is the current minimizer; the upside-down red triangle is the worst agent.  The angles between the two directions and the step sizes are larger for lighter agents. Random descent is a better alternative to gradient descent for some agents and worse for others.}
       \label{fig:direction_comparison}
\end{figure}

Figures \ref{fig:minimizers} and \ref{fig:direction_comparison} provide additional information about the `inner working' of the \SBD dynamics. Figure \ref{fig:minimizers} shows how the \SBD toggles between minimizers and heaviest agents: the loss function decays rapidly for the minimizers. Heavy agents, however, may arise due to merging multiple agents near local minima. The mass of such agents is then slowly transferred to the minimizers with a better minimum.
Figure \ref{fig:direction_comparison} demonstrates the difference between the randomized direction and the gradient direction.

\begin{table}[h]
\setstretch{1.5}
       \centering
    \begin{tabular}{|c| |cc|cc|cc|cc|cc|}
    \hline
    \backslashbox{$d$}{$N$} & \multicolumn{2}{c|}{10} & \multicolumn{2}{c|}{25} & \multicolumn{2}{c|}{50} & \multicolumn{2}{c|}{100} \\
\hline
& \SBD & SBGD & \SBD & SBGD & \SBD & SBGD & \SBD & SBGD\\
\hline
12 & 13.4\% & 24.6\% & 71.1\% & 97.1\% & 100.0\% & 100.0\% & 100.0\% & 100.0\% & \\
14 & \textbf{3.8\%} & 1.3\% & \textbf{60.0\%} & 22.0\% & \textbf{100.0\%} & 49.9\% & \textbf{100.0\%} & 84.1\% & \\
16 & 0.3\% & 0.0\% & \textbf{38.3\%} & 0.1\% & \textbf{95.0\%} & 0.8\% & \textbf{100.0\%} & 1.6\% & \\
18 & 0.0\% & 0.0\% & \textbf{16.3\%} & 0.0\% & \textbf{79.7\%} & 0.0\% & \textbf{99.6\%} & 0.0\% & \\
20 & 0.0\% & 0.0\% & \textbf{1.4\%} & 0.0\% & \textbf{25.1\%} & 0.0\% & \textbf{74.5\%} & 0.0\% & \\
\hline
\end{tabular}
\smallskip
 \caption{Success rate of \SBD vs. SBGD for Ackley test function with mass transfer parameter $\mq=4$.}
 \label{tab:power4}
 \end{table}

\begin{table}[h]
%\vspace*{-0.7cm}
\setstretch{1.5}
       \centering
    \begin{tabular}{|c| |cc|cc|cc|cc|cc|}
    \hline
    \backslashbox{$d$}{$N$} & \multicolumn{2}{c|}{10} & \multicolumn{2}{c|}{25} & \multicolumn{2}{c|}{50} & \multicolumn{2}{c|}{100} \\
\hline
& \SBD & SBGD & \SBD & SBGD & \SBD & SBGD & \SBD & SBGD\\
\hline
12 & 13.0\% & 20.7\% & 75.8\% & 93.3\% & 100.0\% & 99.8\% & 100.0\% & 100.0\% & \\
14 & \textbf{4.2\%} & 1.1\% & \textbf{66.9\%} & 15.3\% & \textbf{100.0\%} & 39.4\% & \textbf{100.0\%} & 77.8\% & \\
16 & 0.1\% & 0.0\% & \textbf{38.4\%} & 0.1\% & \textbf{99.8\%} & 0.8\% & \textbf{100.0\%} & 1.4\% & \\
18 & 0.0\% & 0.0\% & \textbf{14.6\%} & 0.0\% & \textbf{87.3\%} & 0.0\% & \textbf{100.0\%} & 0.0\% & \\
20 & 0.0\% & 0.0\% & \textbf{1.0\%} & 0.0\% & \textbf{30.7\%} & 0.0\% & \textbf{84.7\%} & 0.0\% & \\
\hline
\end{tabular}
\smallskip
 \caption{Success rate of \SBD vs. SBGD for Ackley test function with mass transfer parameter $\mq=8$.}
 \label{tab:power8}
 \end{table}

\begin{table}[h]
%\vspace*{-0.7cm}
\setstretch{1.5}
\centering
\begin{tabular}{|c| |cc|cc|cc|cc|cc|}
\hline
\backslashbox{$d$}{$N$} & \multicolumn{2}{c|}{10} & \multicolumn{2}{c|}{25} & \multicolumn{2}{c|}{50} & \multicolumn{2}{c|}{100} \\
\hline
q  & 8 & 4 & 8 & 4 & 8 & 4 & 8 & 4 \\
\hline
\textbf{Ackley} \\
\hline
14 & \textbf{4.2\%} & 3.8\% & \textbf{66.9\%} & 60.0\% & 100.0\% & 100.0\% & 100.0\% & 100.0\% \\
16 & 0.1\% & 0.3\% & 38.4\% & 38.3\% & \textbf{99.8\%} & 95.0\% & 100.0\% & 100.0\% \\
18 & 0.0\% & 0.0\% & 14.6\% & 16.3\% & \textbf{87.3\%} & 79.7\% & 100.0\% & 99.6\%  \\
20 & 0.0\% & 0.0\% & 1.0\% & 1.4\% & \textbf{30.7\%} & 25.1\% & \textbf{84.7\%} & 74.5\%  \\
\hline
\textbf{Rastrigin} \\
\hline
2 & \textbf{42.5\%} & 37.4\% & 99.0\% & 98.8\% & 100.0\% & 100.0\% & 100.0\% & 100.0\% \\
3 & 6.2\% & 6.5\% & \textbf{32.4\%} & 29.0\% & \textbf{80.1\%} & 74.3\% & \textbf{99.1\%} & 98.0\% \\
4 & 0.8\% & 1.0\% & 4.7\% & 4.9\% & \textbf{14.3\%} & 11.5\% & \textbf{35.2\%} & 30.3\% \\
5 & 0.2\% & 0.1\% & 1.1\% & 0.9\% & \textbf{3.0\%} & 1.7\% & \textbf{4.8\%} & 3.7\% \\
\hline
\textbf{Rosenbrock}\\ 
\hline
%2 & 13.4\% & 13.0\% & 54.0\% & 53.9\% & 95.6\% & 96.3\% & 100.0\% & 100.0\% \\
3 & 2.6\% & 2.8\% & \textbf{12.0\%} & 9.2\% & \textbf{53.7\%} & 45.4\% & \textbf{94.0\%} & 92.0\% \\
4 & 2.2\% & 2.4\% & \textbf{6.6\%} & 5.5\% & \textbf{24.0\%} & 16.7\% & \textbf{63.7\%} & 60.9\% \\
5 & 0.9\% & 1.1\% & 2.4\% & 1.9\% & \textbf{11.8\%} & 7.0\% & \textbf{37.4\%} & 28.5\% \\
6 & 0.5\% & 0.5\% & 1.3\% & 1.1\% & \textbf{6.1\%} & 4.3\% & \textbf{18.1\%} & 13.5\% \\
\hline
\textbf{Styblinski}\\ 
\hline
6 & 9.3\% & 9.3\% & \textbf{42.4\%} & 38.2\% & 73.1\% & 72.9\% & 96.0\% & 95.8\% \\
8 & 2.6\% & 2.9\% & \textbf{12.8\%} & 10.3\% & \textbf{31.2\%} & 30.4\% & \textbf{60.3\%} & 58.4\% \\
10 & 0.8\% & 0.5\% & 4.2\% & 3.5\% & \textbf{11.4\%} & 10.4\% & \textbf{23.5\%} & 20.8\% \\
12 & 0.2\% & 0.1\% & 1.1\% & 1.2\% & 3.6\% & 3.3\% & 7.4\% & 8.1\% \\
\hline
\end{tabular}
\smallskip
 \caption{Success rates of \SBD vs. SBGD for global optimization of the $d$-dimensional objective functions with mass transfer parameter  $\mq=4$ and $\mq=8$. this is to be compared with the corresponding success rate using mass transfer parameter $\mq=2$, reported in the Tables \ref{tab:SBD multiD ackley}, \ref{tab:rastrigin}, \ref{tab:rosenbrock} and \ref{tab:styblinski-tang}.}
 \label{tab:q_powers}
 \end{table}

\subsection{\SBD with higher order mass transition $\mq>2$.}\label{sec:high-q}
In this section we revisit the benchmark examples with different mass transfer parameter $\mq$,
\be\label{eq:mq-transfer}
\left\{\begin{array}{lll}
         \ m_i^{n+1} & =m^n_i  -\etaintwo m^n_i, & i\neq \imin \\ \\
          \ m_\imin^{n+1} & =\displaystyle m_\imin^n +\sum \limits_{i\neq \imin}\etaintwo m_i^{n}, & 
     \end{array} \right. \quad \etain := \Big(\frac{F(\bxin)-F^n_{\textnormal{min}}}{F^{n}_{\textnormal{max}}-F^{n}_{\textnormal{min}}}\Big)^\mq \in (0,1],
\ee
We find that  increasing $\mq$ in the mass transfer protocol
\eqref{eq:mq-transfer}, improves  the success rate of SBRD. Previously, we found $\mq=2$ to be an optimal choice for SBGD. However, as shown in Tables \ref{tab:power4} and \ref{tab:power4} for the Ackley test function, higher $\mq=4$ and respectively $\mq=8$, has a dramatic effect in improving the success rate of \SBD over SBGD. Randomization favors higher transfer parameter $\mq$. 
Indeed, increasing $\mq$ enforces smaller amounts of  mass transfer in \eqref{eq:mq-transfer} so that \SBD becomes more `egalitarian': both the heavier leading agents and the lighter exploring agents are allowed more time (iterations) to settle or to explore, and hence the rate of change  for mass configuration of the swarm become smaller. In particular, this allows a more effective exploration of the random-based descent, improving the overall performance of SBRD. This is demonstrated in Table  \ref{tab:q_powers}.

\pagebreak
\bibliographystyle{plain}
\bibliography{Ref}

\begin{thebibliography}{10}

\bibitem{armijo1966}
Larry Armijo.
\newblock Minimization of functions having lipschitz continuous first partial
  derivatives.
\newblock {\em Pacific Journal of mathematics}, 16(1):1--3, 1966.

\bibitem{boyd2004convex}
Stephen~P Boyd and Lieven Vandenberghe.
\newblock {\em Convex optimization}.
\newblock Cambridge university press, 2004.

\bibitem{CBO-analytical}
Jos{\'e}~A Carrillo, Young-Pil Choi, Claudia Totzeck, and Oliver Tse.
\newblock An analytical framework for consensus-based global optimization
  method.
\newblock {\em Mathematical Models and Methods in Applied Sciences},
  28(06):1037--1066, 2018.

\bibitem{carrillo2021consensus}
Jos{\'e}~A Carrillo, Shi Jin, Lei Li, and Yuhua Zhu.
\newblock A consensus-based global optimization method for high dimensional
  machine learning problems.
\newblock {\em ESAIM: Control, Optimisation and Calculus of Variations}, 27:S5,
  2021.

\bibitem{grassi2023pso}
Sara Grassi, Hui Huang, Lorenzo Pareschi, and Jinniao Qiu.
\newblock Mean-field particle swarm optimization.
\newblock In Benoit~Perthame Weizhu~Bao, Peter A.~Markowich and Eitan Tadmor,
  editors, {\em Modeling and Simulation for Collective Dynamics}, pages
  127--194. World Scientific, 2023.

\bibitem{CBO-semidiscrete}
Seung-Yeal Ha, Shi Jin, and Doheon Kim.
\newblock Convergence of a first-order consensus-based global optimization
  algorithm.
\newblock {\em Mathematical Models and Methods in Applied Sciences},
  30(12):2417--2444, 2020.

\bibitem{CBO-timediscrete}
Seung-Yeal Ha, Shi Jin, and Doheon Kim.
\newblock Convergence and error estimates for time-discrete consensus-based
  optimization algorithms.
\newblock {\em Numerische Mathematik}, 147(2):255--282, 2021.

\bibitem{benchmarks}
Momin Jamil and Xin-She Yang.
\newblock A literature survey of benchmark functions for global optimisation
  problems.
\newblock {\em International Journal of Mathematical Modelling and Numerical
  Optimisation}, 4(2):150--194, 2013.

\bibitem{PSO}
James Kennedy and Russell Eberhart.
\newblock Particle swarm optimization.
\newblock In {\em Proceedings of ICNN'95-international conference on neural
  networks}, volume~4, pages 1942--1948. IEEE, 1995.

\bibitem{kingma2017adam}
Diederik~P Kingma and Jimmy Ba.
\newblock Adam: A method for stochastic optimization.
\newblock {\em arXiv preprint arXiv:1412.6980}, 2017.

\bibitem{SA1}
Scott Kirkpatrick, C~Daniel Gelatt, and Mario~P Vecchi.
\newblock Optimization by simulated annealing.
\newblock {\em science}, 220(4598):671--680, 1983.

\bibitem{liu2022adaptive}
Hailiang Liu and Xuping Tian.
\newblock An adaptive gradient method with energy and momentum.
\newblock {\em Annals of Applied Mathematics}, 38(2):183--222, 2022.

\bibitem{law1965ensembles}
Stanislaw {\L}ojasiewicz.
\newblock Ensembles semi-analytiques.
\newblock {\em IHES notes}, 1965.

\bibitem{lojasiewicz1993geometrie}
Stanislaw {\L}ojasiewicz.
\newblock Sur la g{\'e}om{\'e}trie semi-et sous-analytique.
\newblock In {\em Annales de l'institut Fourier}, volume~43, pages 1575--1595,
  1993.

\bibitem{LTZ2022}
Jingcheng Lu, Eitan Tadmor, and Anil Zenginoglu.
\newblock Swarm-based gradient descent method for non-convex optimization.
\newblock {\em arXiv preprint arXiv:2211.17157}, 2022.

\bibitem{nocedal1999numerical}
Jorge Nocedal and Stephen~J Wright.
\newblock {\em Numerical optimization}.
\newblock Springer, 1999.

\bibitem{nocedal2006conjugate}
Jorge Nocedal and Stephen~J Wright.
\newblock {\em Conjugate gradient methods}.
\newblock Springer, 2006.

\bibitem{palais1964generalized}
Richard~S Palais and Stephen Smale.
\newblock A generalized morse theory.
\newblock {\em Bulletin of the American Mathematical Society}, 70:165--172,
  1964.

\bibitem{CBO1}
Ren{\'e} Pinnau, Claudia Totzeck, Oliver Tse, and Stephan Martin.
\newblock A consensus-based model for global optimization and its mean-field
  limit.
\newblock {\em Mathematical Models and Methods in Applied Sciences},
  27(01):183--204, 2017.

\bibitem{reynolds1987flocks}
Craig~W Reynolds.
\newblock Flocks, herds and schools: A distributed behavioral model.
\newblock In {\em Proceedings of the 14th annual conference on Computer
  graphics and interactive techniques}, pages 25--34, 1987.

\bibitem{tadmor2021mathematics}
Eitan Tadmor.
\newblock On the mathematics of swarming: emergent behavior in alignment
  dynamics.
\newblock {\em Notices of the AMS}, 68(4):493--503, 2021.

\bibitem{totzeck2022trends}
Claudia Totzeck.
\newblock Trends in consensus-based optimization.
\newblock In N.~Bellomo J.~A. Carrillo and E.~Tadmor, editors, {\em Active
  Particles, Volume 3}, pages 201--226. Springer, 2022.

\bibitem{SA2}
Peter~JM Van~Laarhoven and Emile~HL Aarts.
\newblock Simulated annealing.
\newblock In {\em Simulated annealing: Theory and applications}, pages 7--15.
  Springer, 1987.

\bibitem{wolfe1969}
Philip Wolfe.
\newblock Convergence conditions for ascent methods.
\newblock {\em SIAM review}, 11(2):226--235, 1969.

\end{thebibliography}

\end{document}